\documentclass[preprint, 11pt]{elsarticle}

\usepackage[T1]{fontenc}
\usepackage[utf8]{inputenc}
\usepackage[hmargin=0.9in]{geometry}
\usepackage[babel,german=quotes]{csquotes}
\usepackage{subfig}
\usepackage{graphicx}
\usepackage[colorlinks=true,citecolor=blue,urlcolor=blue]{hyperref}
\usepackage{caption}
\usepackage{amsmath}          
\usepackage{amsthm}           
\usepackage{amssymb}          
\usepackage{bm, amsfonts}     
\usepackage{xcolor}
\usepackage{fixmath}
\usepackage{tikz}
\usepackage{bm}
\usepackage{makecell}

\graphicspath{{./pictures/small/}}

\newdefinition{cor}{Corollary}
\newdefinition{example}{Example}
\newproof{pf}{Proof}

\newcommand {\dx} {\,{\rm d}{\mathbf x}}

\newcommand {\ds} {\,{\rm d}{\mathrm s}}

  \newcommand{\R}{\mathbb{R}}
  \newdefinition{rmk}{Remark}
  
  \newcommand{\pd}[2]{\frac{\partial #1}{\partial #2}}

  \newcommand{\td}[2]{\frac{\mathrm d #1}{\mathrm d #2}}

\newcommand{\beq}{\begin{equation}}
\newcommand{\eeq}{\end{equation}}

\newcommand{\Heps}{H_\epsilon}
\newcommand{\Seps}{S_\epsilon}
\newcommand{\bfv}{\mathbf{v}}

\DeclareMathOperator\erf{erf}

\newcommand\red[1]{#1} 
\newcommand\blue[1]{#1}	
\newcommand\cyan[1]{#1}  

\newcommand{\bfq}{{\bf q}}

\makeatletter
\def\ps@pprintTitle{%
  \let\@oddhead\@empty
  \let\@evenhead\@empty
  \def\@oddfoot{
    \footnotesize\itshape
    \hfill\today
  }%
  \let\@evenfoot\@oddfoot}
\makeatother

\begin{document}

\begin{frontmatter} 
  \title{An unfitted finite element method using level set functions for\\
    extrapolation into deformable diffuse interfaces}

\author{Dmitri Kuzmin\corref{cor1}}
\ead{kuzmin@math.uni-dortmund.de}

\cortext[cor1]{Corresponding author}

\author{Jan-Phillip B\"acker}
\ead{jan-phillip.baecker@tu-dortmund.de}

\address{Institute of Applied Mathematics (LS III), TU Dortmund University\\ Vogelpothsweg 87,
  D-44227 Dortmund, Germany}

\journal{Journal of Computational Physics}

\begin{abstract}  

  We explore a new way to handle flux boundary conditions imposed on level sets. The proposed approach is a diffuse interface version of the shifted boundary method (SBM) for continuous Galerkin discretizations of conservation laws in embedded domains. We impose the interface conditions weakly and approximate surface integrals by volume integrals. The discretized weak form of the governing equation has the structure of an immersed boundary finite element method. \red{That is, integration is performed over a fixed fictitious domain. Source terms are included to account for interface conditions and extend the boundary data into the complement of the embedded domain. The calculation of these extra terms} requires (i) construction of an approximate delta function and (ii) extrapolation of embedded boundary data into quadrature points. We accomplish these tasks using a level set function, which is given analytically or evolved numerically. A~globally defined averaged gradient of this approximate signed distance function is used to construct a simple map to the closest point on the interface. The normal and tangential derivatives of the numerical solution at that point are calculated using the interface conditions and/or interpolation on uniform stencils. Similarly to SBM, extrapolation of data back to the quadrature points is performed using Taylor expansions. Computations that require extrapolation are restricted to a narrow band around the interface. Numerical results are presented for elliptic, parabolic, and hyperbolic test problems, which are specifically designed to assess the error caused by the numerical treatment of interface conditions on fixed and moving boundaries in 2D.
  
\end{abstract}
\end{frontmatter}

\pagebreak

\section{Introduction}
\label{sec:intro}

Many unfitted mesh (fictitious domain, immersed boundary, cut cell) discretization methods have been proposed in the literature to enforce interface conditions for applications ranging from the Poisson equation to multiphase flow models and fluid-structure interaction problems. In this work, we focus on weak formulations in which interfacial terms are represented by surface or volume integrals. In level set methods \cite{osher,sethian_review}, the interface is defined as a manifold on which an approximate signed distance function (SDF) becomes equal to zero. The SDF property of advected level sets is usually preserved using a postprocessing procedure known as \emph{redistancing}. The monolithic conservative version proposed in \cite{mono} incorporates redistancing and mass correction into a nonlinear transport equation for the level set function. A variety of unfitted finite element methods (FEM) are available for numerical integration over sharp, diffuse, or surrogate interfaces. A sharp interface method enriches the local finite element space to accommodate boundary/jump conditions on an embedded surface. For example, this strategy is used in XFEM \cite{xfem}, unfitted Nitsche FEM \cite{hansbo2,wadbro}, and CutFEM \cite{cutfem,hansbo2014} approaches.~Integration over sharp interfaces has the potential drawback that small cut cells may cause ill-conditioning or severe time step restrictions. Possible remedies to this problem include the use of stabilization terms \cite{hansbo2014,zahedi}, ghost penalties \cite{burman,olsh}, and (explicit-)implicit schemes for time integration \cite{may}.

Diffuse interface methods replace surface integrals by volume integrals depending on regularized Dirac delta functions \cite{kublik2016,kublik2018,teigen2011,teigen2009}. This family of unfitted FEM is closely related to classical immersed boundary methods \cite{mittal,peskin2002} and well suited, e.g., for implementation of surface tension effects in two-phase flow models \cite{CSF,hysing}. The construction of approximate delta functions for level set methods was addressed in \cite{engquist,smereka,towers,zahedi}. M\"uller et al.~\cite{mueller2013} devised an elegant alternative, which uses the divergence theorem and divergence-free basis functions to reduce integration over an embedded interface to that over a fitted boundary. Another promising new approach is the \emph{shifted boundary method} (SBM) \cite{sbm,sbm1,sbm2}, which belongs to the family of surrogate interface approximations. To avoid the need for dealing with small cut cells, the SBM applies surrogate boundary conditions on common edges/faces of cut and uncut cells. The extrapolation of data to shifted boundaries is performed using Taylor expansions. Interface jump conditions are treated similarly \cite{sbm}. The SBM approach is conceptually simple and \cyan{justified} by theoretical analysis \cite{atallah,sbm1}. A~potential drawback is high complexity of closest-point projection algorithms for constructing a map to the true interface.

In the present paper, we introduce a diffuse interface version of SBM, in which surface integration is performed using an algorithm that fits into the framework developed in \cite{kublik2016}. Approximating surface integrals by volume integrals, we use level set functions to construct regularized delta functions, approximate global normals, and closest-point mappings. The extrapolation of interface data into the quadrature points for volume integration involves the following steps: (i) find the closest point on the interface; (ii) calculate the quantities of interest at that point; (iii) perform constant extension of normal derivatives (for diffusive fluxes) and linear extension of solution values (for inviscid fluxes). In Step~(i), we evaluate the approximate normal at the quadrature point, construct a straight line parallel to the resulting vector, and perform a line search to find the nearest root of the level set function. In Step~(ii), normal derivatives are calculated using the given interface value and the interpolated value at the internal edge of the diffuse interface. If tangential derivatives are needed for linear extrapolation in Step~(iii), they are approximated using an interpolation stencil centered at the internal point. The way in which the interface data is projected into quadrature points distinguishes our algorithm from shifted boundary methods and extrapolation-based approaches proposed in \cite{kublik2016,kublik2018}.

The design of our diffuse level set method is motivated by embedded domain problems in which the interface motion is driven by gradients of concentration fields \cite{hogea,chen,teigen2011,teigen2009,vermolen2007}. Level set algorithms for such problems require a suitable extension of the interface velocity \cite{sethian,chen,utz}. We perform it in the same way as extrapolation of normal fluxes. The level set function is evolved as in \cite{mono,monons}. The support of extended interface data is restricted to a narrow band around the interface. The cost of integrating extrapolated quantities is minimized in this way. For numerical experiments, we design a set of two-dimensional elliptic, parabolic, and hyperbolic test problems with closed-form analytical solutions. The interface is fixed in the elliptic and hyperbolic test cases but is moving in the parabolic scenario. We report and discuss the results of grid convergence studies for our new benchmarks. 

\red{
\section{Interface description}

Let $\Omega_+(t)\subset\R^d,\ d\in\{1,2,3\}$ be a time-dependent domain which is enclosed by an evolving interface $\Gamma(t)$. In our method, $\Omega_+(t)$ remains fully embedded into a fixed fictitious domain $\Omega\subset\R^d$ for all times $t\ge 0$. The instantaneous position of 
$\Gamma(t)$ is determined by a level set function $\phi(\mathbf{x},t)$ which
 is positive for $\mathbf{x}\in\Omega_+(t)$ and negative in
$\Omega_-(t)=
\bar\Omega\backslash\bar\Omega_+(t)$. This assumption implies that
\beq
\Gamma(t)=\{\mathbf x\in \Omega\,:\, \phi(\mathbf x,t)=0\}.
\eeq
The vector fields $\mathbf{n}_\pm=\mp
\frac{\nabla\phi}{|\nabla\phi|}$ represent extended unit outward
normals to $\partial\Omega_\pm(t)\cap
\Gamma(t)$. If $\phi$ is a signed distance function satisfying
the Eikonal equation $|\nabla\phi|=1$, then $\mathbf{n}_-
=\nabla\phi=-\mathbf{n}_+$.

Classical level set methods \cite{osher,sethian_review} evolve
$\phi$ by solving the linear advection equation
\beq\label{ls:sharp}
\pd{\phi}{t}+\mathbf{v}\cdot\nabla\phi=0 \qquad\mbox{in}\ \Omega.
\eeq
In incompressible two-phase flow models, the velocity $\mathbf{v}$
is obtained by solving the Navier-Stokes equations
with piecewise-constant density and viscosity \cite{hysing,monons}. In mass
and heat transfer models,
$\mathbf{v}$ may represent an extension of a normal velocity depending
on concentration gradients \cite{hogea,vermolen2007}. 
The exact solution of \eqref{ls:sharp}
satisfies the  nonlinear transport equation
\beq\label{vof:noncons}
\pd{H(\phi)}{t}
+\mathbf{v}\cdot \nabla H(\phi)=0 \qquad\mbox{in}\ \Omega,
\eeq
where the derivatives of the discontinuous Heaviside function
\beq\label{heaviside}
H(\phi)=\begin{cases}
1 & \text{ if } \phi > 0,\\
0 & \text{ if } \phi < 0
\end{cases}
\eeq
should be understood in the sense of distributions. If the velocity
field $\mathbf{v}$ is divergence-free, then 
\beq\label{vof:sharp}
\pd{H(\phi)}{t}
+ \nabla\cdot(\mathbf{v}H(\phi))=0 \qquad\mbox{in}\ \Omega.
\eeq
Integrating \eqref{vof:sharp} over $\Omega$, using the
divergence theorem and the fact that $H(\phi)|_{\partial\Omega}=0$
 for fully embedded domains $\Omega_+(t)$,
we find that the volume $|\Omega_+(t)|=\int_\Omega H(\phi(\mathbf{x},t))\dx$ of the
embedded domain is independent of $t$ if $\nabla\cdot\mathbf{v}=0$ 
in $\Omega$. However, approximate solutions of \eqref{ls:sharp}
may fail to satisfy a discrete form of the local
conservation law for $H(\phi)$. 
This drawback of the level set approach can be cured using various
mass correction procedures. The monolithic level set method that we
present in Section~\ref{sec:diffuse} is conservative by construction and does not
require any postprocessing (such as redistancing, which is commonly
employed to approximately preserve the distance function property).

\section{Sharp interface problem}
\label{sec:sharp}

As a representative model of practical interest, we first consider the parabolic conservation law
\beq\label{pde}
\pd{u}{t}+\nabla\cdot [\mathbf{f}(u)-\kappa\nabla u]=0\qquad \mbox{in}\ \Omega_+(t),
\eeq
where $\mathbf{f}(u)$ is an inviscid flux, $\kappa\ge 0$ is a diffusion
coefficient, and $t>0$ is a fixed instant.
The boundary and initial conditions that we impose weakly in this work are given by
 \begin{align}
 [\mathbf{f}(u)-\kappa\nabla u]\cdot\mathbf{n}_+
=g_\Gamma(u,u_\Gamma) &\quad\mbox{on}\ \Gamma(t),\label{bc}\\
 u=u_0  &\quad\mbox{in}\ \Omega_+(0),\label{ic}
 \end{align}
 where $u_0$ is the initial data,
 $u_\Gamma$ is the boundary data and for $\mathbf{x}\in\Gamma(t)$
 \begin{align}
   g_\Gamma(u,u_\Gamma)(\mathbf{x})=
   f_\mathbf{n}(u(\mathbf{x}),u_\Gamma(\mathbf{x}))
   -\kappa\lim_{\epsilon\searrow 0}
   \frac{u_\Gamma(\mathbf{x})-u(\mathbf{x}-
     \epsilon\mathbf{n}_+(\mathbf{x}))}{\epsilon}
 \end{align}
is the total flux in the direction of the unit outward
 normal $\mathbf{n}_+$. The component
 $f_\mathbf{n}(u(\mathbf{x}),u_\Gamma(\mathbf{x}))$ is
 defined by a Riemann solver. The normal derivative
 is written as the limit of a ratio depending on $u_\Gamma(\mathbf{x})$
 because we will approximate this limit by a finite difference in 
 Section \ref{sec:fluxes}.

 To formulate our sharp interface problem in
 the fictitious domain $\Omega$, we also need
 the  condition
 \beq\label{ghost}
u=u_-\qquad \mbox{in}\ \Omega_-(t),
\eeq
where $u_-$ is a suitable extension of $u_\Gamma$, such as the
Dirichlet-Neumann \emph{harmonic extension} \cite{gorb}
\begin{alignat*}{2}
-\Delta u_-&=0\ &&\quad\mbox{in}
  \ \Omega_-(t),\\ u_-&=u_\Gamma\ &&\quad\mbox{on}\ \Gamma(t),\\
  \mathbf{n}\cdot
  \nabla u_-&=0 && \quad\mbox{on}\ \partial\Omega.
\end{alignat*}
In the Neumann boundary condition,
$\mathbf{n}$ denotes the unit outward normal to
the boundary $\partial\Omega$ of $\Omega$. Another way to
define $u_-$ will be discussed in Section \ref{sec:ghost}.
The existence of continuous extension operators from a bounded
 domain  to $\R^d$ is guaranteed by the theoretical result
 referred to in \cite{olsh}. 

For illustration purposes, let us first notice that for 
smooth $u=u(\mathbf{x},t)$ and a smooth
test function $w=w(\mathbf{x})$, equations \eqref{pde}--\eqref{ghost}
of the strong form fictitious domain problem imply
\begin{align}\label{weak-fd}
 \int_{\Omega_+(t)}w\pd{u}{t}\dx
  &- \int_{\Omega_+(t)}\nabla w\cdot
     [\mathbf{f}(u)-\kappa\nabla u]\dx
     +\int_{\Omega_-(t)}\gamma_\Omega w(u-u_-)\dx\nonumber\\
      &
      + \int_{\Gamma(t)}w g(u,u_\Gamma)\ds=0
\end{align}
for any positive $\gamma_\Omega$ and fixed $t>0$. The integral
containing $u_-$ is a \emph{ghost penalty} term, It represents
a weighted residual of \eqref{ghost}. The other terms
represent a weighted residual of \eqref{pde} after
integration by parts and substitution of the flux boundary
conditions \eqref{bc}.
Since $w$ is independent of $t$, the application of the Reynolds
transport theorem  to the first integral on the left-hand side
of \eqref {weak-fd} yields
\beq\label{rtt}
\int_{\Omega_+(t)}w\pd{u}{t}\dx=
\int_{\Omega_+(t)}\pd{(wu)}{t}\dx=\td{}{t}
\int_{\Omega_+(t)}wu\dx-\int_{\Gamma(t)}wuv_n\ds,
\eeq
where $v_n=\mathbf{v}\cdot\mathbf{n}_+$ is the velocity of the
interface in the normal direction $\mathbf{n}_+$.
The assumption that $u$ and $w$ are smooth will be weakened at
the level of the finite element discretization.

Using the Heaviside function defined by \eqref{heaviside}, volume integration over $\Omega_\pm(t)$ can be replaced by that over the fixed
fictitious domain $\Omega$, Substituting \eqref{rtt} into \eqref{weak-fd}, we obtain
\begin{align}
 \td{}{t} \int_{\Omega}H(\phi)wu\dx
  &- \int_{\Omega}H(\phi)\nabla w\cdot
    [\mathbf{f}(u)-\kappa\nabla u]\dx+\int_{\Omega}\gamma_\Omega w(u-u_\Omega)\dx
    \nonumber\\
      &
      + \int_{\Gamma(t)}w[g(u,u_\Gamma)-v_nu]\ds
=0,\label{weak2}
\end{align}
where $u_\Omega=H(\phi)u+(1-H(\phi))u_-$.
 For our choice of $u_\Omega$, the contribution of
$\int_\Omega\gamma_\Omega(u-u_\Omega)\dx$
        imposes the Dirichlet extension condition
        \eqref{ghost} weakly in the external subdomain
        $\Omega_-(t)$, where $u_\Omega=u_-$, and
        vanishes in $\Omega_+(t)$, where
        $u_\Omega=u$. Discrete counterparts of such ghost penalties are used
in unfitted finite element methods for stabilization and regularization purposes
\cite{burman,olsh}.

 \begin{rmk}\label{rem:ghost}
Neumann-type ghost penalties of the form
$\int_{\Omega}\gamma_\Omega \nabla w\cdot(\nabla u-\mathbf{g}_\Omega)\dx$
can be defined using a vector field $\mathbf{g}_\Omega$ such that
$\mathbf{g}_\Omega=\nabla u$ in $\Omega_+(t)$. This version penalizes
the weak residual of
\begin{alignat*}{2}
-\Delta u&=-\nabla\cdot \mathbf{g}_\Omega\ &&\quad\mbox{in}
  \ \Omega_-(t),\\ u&=u_\Gamma\ &&\quad\mbox{on}\ \Gamma(t),\\
  \mathbf{n}\cdot
  \nabla u_-&=0 && \quad\mbox{on}\ \partial\Omega.
\end{alignat*}
If $\nabla\cdot\mathbf{g}_\Omega=0$  in
$\Omega_-(t)$, a harmonic extension
of $u_\Gamma$ into $\Omega_-(t)$
is defined by the solution of
 this boundary value problem.
 In Section \ref{sec:ghost}, we
 construct a discrete counterpart of $\mathbf{g}_\Omega$
 using extrapolation. Other kinds of discrete ghost penalties
 for fictitious domain methods can be found in \cite{olsh}.
 \end{rmk}
 }
 
\section{Diffuse interface problem}
\label{sec:diffuse}

To avoid numerical difficulties associated with volume integration of
discontinuous functions and surface integration over sharp evolving interfaces,
we approximate \eqref{weak2} by (cf. \cite{hysing,kublik2016,zahedi})
\begin{align}
 \td{}{t} \int_{\Omega}H_\epsilon(\phi)wu\dx
  &- \int_{\Omega}H_\epsilon(\phi)\nabla w\cdot
    [\mathbf{f}(u)-\kappa\nabla u]\dx+\int_{\Omega}\gamma_\Omega w(u-u_\Omega)\dx
    \nonumber\\
      &
    + \int_{\Omega}
    wG(\phi,u,u_\Gamma)
        \delta_\epsilon(\phi)|\nabla\phi|\dx
=0\qquad \forall w\in V(\Omega),\label{weak3}
\end{align}
where $G(\phi,u,u_\Gamma)$ is an extension of the interface flux
\cyan{$g(u,u_\Gamma)-v_nu$} into $\Omega$.
The functions
$H_\epsilon(\phi)$ and $\delta_\epsilon(\phi)$ represent
regularized approximations to 
$H(\phi)$ and the Dirac delta function
$\delta_\Gamma$, respectively.
The design and analysis of such approximations have received significant
attention in the literature during the last two decades
\cite{engquist,kublik2016,kublik2018,smereka,towers,zahedi}. Many diffuse interface methods
\cite{teigen2011,teigen2009} and level set algorithms \cite{hysing,mono,monons}
use regularized Heaviside and/or delta functions.
The novelty of our unfitted finite element method lies in the construction of
$G(\phi,u,u_\Gamma)$, which we discuss in the next section.
\smallskip

In the numerical experiments of Section~\ref{sec:num}, we use the
 Gaussian regularization \cite[eq. (52)]{zahedi} 
 \beq\label{delta:gauss}
 H_\epsilon(\phi) = \frac{1}{2}\left(1+\erf\left(\frac{\pi \phi}{3\epsilon}\right)\right),
  \qquad \delta_\epsilon(\phi) =\cyan{H_\epsilon'(\phi)}=
  \frac{1}{\epsilon}\sqrt{\frac{\pi}{9}}\exp\left(\frac{-\pi^2\phi^2}{9\epsilon^2}\right),
  \eeq
  where $\erf(x)=\frac{2}{\sqrt{\pi}}\int_0^x e^{-y^2}\mathrm{d}y$ is the
  Gauss error function.

To evolve $\phi$ in moving interface scenarios, we use \cyan{an extension of the
monolithic conservative level set method 
\cite{mono} to general velocity fields $\mathbf{v}$.
The nonlinear evolution
equation for $\phi$ becomes}
\beq\label{phasefield}
\pd{ \Seps(\phi)}{t}+\cyan{\mathbf{v}\cdot\nabla\Seps(\phi)-
\lambda\nabla\cdot
\left(\nabla\phi
-\mathbf{q}\right)}=0\qquad\mbox{in}\ \Omega,
\eeq
where $\Seps(\phi)= 2\Heps(\phi)-1$ is a smoothed
sign function and
$\mathbf{q}=\frac{\nabla\phi}{|\nabla\phi|}$
is the Eikonal flux.
The term depending on a
penalty parameter $\lambda>0$ regularizes
\eqref{phasefield}  and enforces the
approximate distance function property
of $\phi$. The
finite element method for solving
\eqref{phasefield} uses the mixed weak form \cite{mono}
\begin{align}
\int_\Omega w\pd{ \Seps(\phi)}{t}\dx
&+\int_\Omega 
\cyan{\left[w\mathbf{v}\cdot\nabla\Seps(\phi)
  +\lambda\nabla w\cdot(\nabla\phi-\bfq)\right]}\dx
= 0,\label{weakS}
\end{align}
\beq\label{weakq}
\int_\Omega w	\sqrt{|\nabla\phi|^2+\sigma^2}\,\bfq\dx
=\int_\Omega w \nabla\phi\dx.
\eeq
\cyan{Note that the integrand $w\mathbf{v}\cdot\nabla\Seps(\phi)$ vanishes
in regions where $\Seps(\phi)$ is constant. Therefore, it is sufficient
to define the extended velocity $\mathbf{v}$ of the interface in a
narrow band around $\Gamma(t)$.}
If the level
set function $\phi$ satisfies
the Eikonal equation $|\nabla\phi|=1$, then $\mathbf{q}=\nabla\phi$ for
$\sigma=0$. A~nonvanishing small value of the regularization
parameter $\sigma$ is used in \eqref{weakq}
to avoid ill-posedness at critical points, i.e., in the limit $|\nabla\phi|\to 0$.
 For a detailed presentation,
we refer the interested reader to \cite{mono,monons}.

\section{Finite element discretization}

We discretize \eqref{weak3}, \eqref{weakS}, and \eqref{weakq} in space using
the continuous Galerkin method and a conforming mesh $\mathcal T_h=\{
K^{1},\ldots,K^{E_h}\}$ of linear
($\mathbb{P}_1$) or multilinear ($\mathbb{Q}_1$) Lagrange
finite elements. The corresponding global basis functions
are denoted by $\varphi_1,\ldots,\varphi_{N_h}$. They span
a finite-dimensional space \cyan{$V_h$} and
have the property that $\varphi_j(\mathbf{x}_i)=\delta_{ij}$ for
each vertex $\mathbf{x}_i$ of $\mathcal T_h$.
The indices of elements that meet at $\mathbf{x}_i$ are stored in
the integer set $\mathcal E_i$. The global indices
of nodes belonging to $K^e\in\mathcal T_h$ are stored
in the integer set $\mathcal N^e$. The computational
stencil of node $i$ is the index set
$\mathcal N_i=\bigcup_{e\in\mathcal E_i}\mathcal N^e$.

The finite element approximations $(u_h,\phi_h,\mathbf{q}_h)
\in V_h\times V_h\times (V_h)^d$ to $(u,\phi,\mathbf{q})$
can be written as
$$
u_h(\mathbf{x},t)=\sum_{j=1}^{N_h}u_j(t)\varphi_j(\mathbf{x}),\qquad
\phi_h(\mathbf{x},t)=\sum_{j=1}^{N_h}\phi_j(t)\varphi_j(\mathbf{x}),\qquad
\mathbf{q}_h(\mathbf{x},t)=\sum_{j=1}^{N_h}\bfq_j(t)\varphi_j(\mathbf{x}).
$$
The standard Galerkin discretization of our system yields the semi-discrete problem
\begin{align}
 \td{}{t} \int_{\Omega_h}H_\epsilon(\phi_h)w_hu_h\dx
  &- \int_{\Omega_h}H_\epsilon(\phi_h)\nabla w_h\cdot
    [\mathbf{f}(u_h)-\kappa\nabla u_h]\dx
    +\int_{\Omega_h}\gamma_{\Omega,h} w_h(u_h-u_{\Omega,h})\dx
    \nonumber\\
      &
    + \int_{\Omega_h}
    w_hG(\phi_h,u_h,u_\Gamma)\delta_\epsilon(\phi_h)|\nabla\phi_h|\dx
=0\qquad \forall w_h\in V_h,\label{weakh}
\end{align}
\begin{align}
  \int_{\Omega_h}w_h\pd{\Seps(\phi_h)}{t}\dx&+\int_{\Omega_h}
  \cyan{\left[
    w_h\bfv_h\cdot\nabla\Seps(\phi_h)+\lambda_h
 \nabla w_h\cdot(\nabla\phi_h-\bfq_h)\right]}\dx
= 0\qquad \forall w_h\in V_h,\label{weakSh}
\end{align}
\beq\label{weakqh}
\int_{\Omega_h} w_h	\sqrt{|\nabla\phi_h|^2+\sigma^2}\,\bfq_h\dx
=\int_{\Omega_h} w_h \nabla\phi_h\dx\qquad \forall w_h\in V_h,
\eeq
where $\Omega_h=\bigcup_{e=1}^{E_h}K^e$. The diffuse interface thickness is
given by $\epsilon=mh,\ m\ge 1$. Thus the sharp interface
formulation is recovered
in the limit $h\to 0$. The penalty functions $\gamma_{\Omega,h}$ and
$\lambda_{h}$ are chosen to be
piecewise constant. Their values in $K^e$
depend on the local mesh size $h^e$.
If necessary, the semi-discrete problem for $u_h$ can be
stabilized, e.g., using algebraic flux correction tools \cite{afc1}.
No stabilization is needed for
\eqref{weakSh} because $\Seps(\phi_h)$ attains values in $[-1,1]$ by
definition. The coefficients $\bfq_j=(q_j^{(1)},\ldots,q_j^{(d)})$ of $\bfq_h$
 can be approximated using the explicit lumped-mass formula \cite[eq. (20)]{mono}
$$
 q_j^{(k)}=\frac{\int_{\Omega_h}\pd{\phi_h}{x_k}\varphi_j\dx}{
   \int_{\Omega_h}\sqrt{|\nabla\phi_h|^2+\sigma^2}\varphi_j\dx},
 \qquad j=1,\ldots,N_h,\quad
 k=1,\ldots,d.
 $$
 Further implementation details and parameter settings for \eqref{weakSh},\eqref{weakqh}
 can be found in
 \cite{mono,monons}. In Sections \ref{sec:closestpoint}--\ref{sec:ghost}, we present
 the algorithms for calculating $G(\phi_h,u_h,u_\Gamma)$ and $u_{\Omega,h}$.
 The calculation of an extension velocity $\bfv_h$ for the discretized phase
 field equation \eqref{weakSh} is discussed in Section \ref{sec:extvel}.

 \begin{rmk}
 An appropriate choice of a time integrator for the spatial semi-discretization \eqref{weakh} is problem
 dependent. In Section~\ref{sec:num}, we use Heun's
 explicit Runge--Kutta method
 in the hyperbolic test and the implicit Crank--Nicolson scheme in the
 parabolic one. In the elliptic test, we
 march the solution of the corresponding parabolic problem 
 to a steady state using a pseudo-time stepping method.
\end{rmk}

\section{Extrapolation using level sets}
\label{sec:extrap}

The main highlight of our diffuse interface method is a simple
new algorithm for calculating the fluxes $G(\phi_h,u_h,u_\Gamma)$,
ghost penalty functions $u_{\Omega,h}$, and extension velocities
$\mathbf{v}_h$. We define
\beq\label{Gflux}
G(\phi_h,u_h,u_\Gamma)=\cyan{F-\kappa\partial_nU-VU},\qquad
\mathbf{v}_h=-V\mathbf{q}_h
\eeq
using continuous extensions $U,F,\partial_nU,V\in C(\bar\Omega_h)$
of pointwise solution values, inviscid fluxes, normal derivatives, and normal
velocities defined on the zero level set
$\Gamma_h(t)=\{\mathbf{x}\in\bar\Omega_h\,:\, \phi_h(\mathbf{x},t)=0\}$.
Of course, the extensions also depend on the mesh size but
we omit the subscript $h$ for brevity.

\subsection{Closest-point search}
\label{sec:closestpoint}

Let $\phi_h$ be an approximate signed distance function and
$\mathbf{q}_h$ an approximation to the
extended normal $\mathbf{n}_-=\frac{\nabla\phi}{|\nabla\phi|}
=-\mathbf{n}_+$. In our
diffuse level set method, $\phi_h$ and $\mathbf{q}_h$ are
given by \eqref{weakSh} and \eqref{weakqh}, respectively.
However, any other approximation or analytical formula
may be used instead. To calculate the extensions
$U,\partial_nU,$ and/or
$V$ for a quadrature point $\mathbf{x}_Q\in\bar\Omega_h$,
we first need to find the closest 
$\mathbf{x}_\Gamma\in\Gamma_h(t)$ that is connected
to  $\mathbf{x}_Q$ by a line parallel to
$\mathbf{n}_Q:=-\mathbf{q}_h(\mathbf{x}_Q)
\approx\mathbf{n}_+(\mathbf{x}_Q)$. The point
$$\mathbf{x}_\Gamma^*=\mathbf{x}_Q
+\phi_h(\mathbf{x}_Q)\mathbf{n}_Q$$
is usually a good approximation. Indeed, if $\phi_h$ is an
exact SDF, then $\phi_h(\mathbf{x}_\Gamma^*)=0$ and, therefore,
$\mathbf{x}_\Gamma=\mathbf{x}_\Gamma^*$~is the closest point.
If $\phi_h$ is not exact, $\mathbf{x}_\Gamma$ can be found
using a line search along
$$\hat{\mathbf{x}}(\xi)=\mathbf{x}_Q+\xi\,
\mathrm{sign}(\phi_h(\mathbf{x}_Q))
\mathbf{n}_Q,\qquad\xi\in\R.$$
Note that $\hat{\mathbf{x}}(0)=\mathbf{x}_Q$ and
$\hat{\mathbf{x}}(\xi_\Gamma^*)=\mathbf{x}_\Gamma^*$
for $\xi_\Gamma^*=|\phi_h(\mathbf{x}_Q)|$. Let
$\hat\phi_h(\xi)=\phi_h(\hat{\mathbf{x}}(\xi))$ for
$\xi\in\R$.
Using this parametrization, the exact location of the
interface point
 $\mathbf{x}_\Gamma$ can be found as follows:

\begin{itemize}
\item Choose the smallest $m\in\mathbb{N}$
   such that $\hat \phi_h(\xi_\Gamma^*+mh)\hat\phi_h(\xi_\Gamma^*-mh)<0$.
 \item Use the bisection method to find the root $\xi_\Gamma\in
[\xi_\Gamma^*-mh,\xi_\Gamma^*+mh]$ of $\hat\phi_h(\xi)$.
\item Set $\mathbf{x}_\Gamma=\hat{\mathbf{x}}(\xi_\Gamma)$.
\end{itemize}
This way to zoom in on $\mathbf{x}_\Gamma$ requires repeated evaluation
of the level set function $\phi_h(\mathbf{x})$ at random locations
on the straight line $\hat{\mathbf{x}}(\xi)$.
The  mesh cells containing the interpolation points are easy
to find on uniform meshes. However, the cost of searching 
becomes large if the mesh is unstructured. Moreover, the
bisection method may require many iterations to achieve
the desired accuracy.
\smallskip

To locate the point $\mathbf{x}_\Gamma$ efficiently on general
meshes, we use the following search algorithm:
\begin{itemize}
\item Set $\xi_0=0$.
\item For $i=1,2,\ldots$ let $\hat{\mathbf{x}}(\xi_i)$ with
  $\xi_i>\xi_{i-1}$ be 
  the next intersection of $\hat{\mathbf{x}}(\xi)$ with the
  boundary $\partial K$ of a mesh cell $K\in\mathcal T_h$. Exit the loop when 
  $\hat \phi_h(\xi_{i-1})\hat\phi_h(\xi_i)$
  becomes negative for $i=m$.
 \item Solve a linear/quadratic equation to find the root $\xi_\Gamma\in
[\xi_{m-1},\xi_m]$ of $\hat\phi_h(\xi)$.
\item Set $\mathbf{x}_\Gamma=\hat{\mathbf{x}}(\xi_\Gamma)$.
\end{itemize}

\begin{figure}[h!]
\centering
\begin{minipage}[t]{0.45\textwidth}
  \centering (a) closest-point search 
\includegraphics[width=\textwidth]{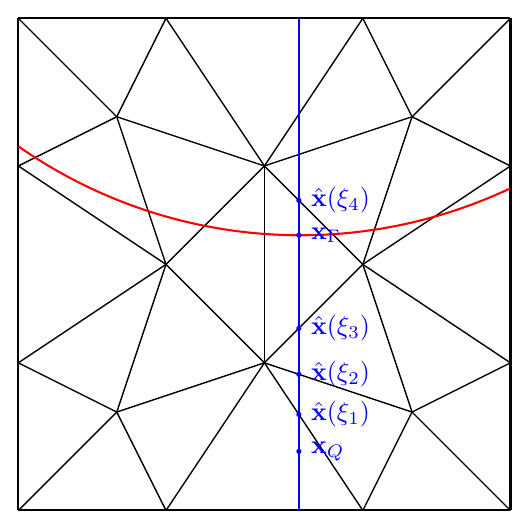}
\end{minipage}%
\begin{minipage}[t]{0.45\textwidth}
  \centering (b) interpolation stencil
  \includegraphics[width=\textwidth]{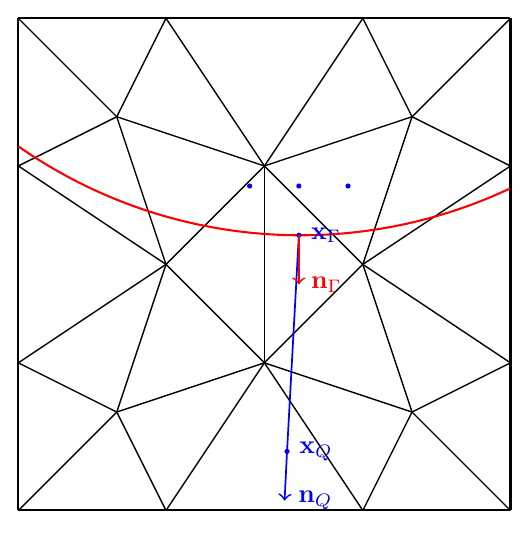}
\end{minipage}
\vskip0.25cm 
\caption{Extrapolation into $\mathbf{x}_Q$ on a triangular mesh.}
\label{sketch}
\end{figure}

 The sketch in
Fig.~\ref{sketch}(a) illustrates the construction of
$\{\hat{\mathbf{x}}(\xi_i)\}_{i=0}^m$ for a  point
$\mathbf{x}_Q=\hat{\mathbf{x}}(0)$ on a triangular mesh.
For $\hat{\mathbf{x}}(\xi_{i-1})$
belonging to $K\in\mathcal T_h$, it is easy to find 
$\hat{\mathbf{x}}(\xi_{i})\in\partial K$ and an adjacent
cell $K'\in\mathcal T_h$ into which the interface navigation
vector $\mathbf{p}_Q=\mathrm{sign}(\phi_h(\mathbf{x}_Q))\mathbf{n}_Q$
is pointing at $\hat{\mathbf{x}}(\xi_{i})$. By definition,
the points $\hat{\mathbf{x}}(\xi_{m-1})$ and
$\hat{\mathbf{x}}(\xi_m)$ belong to the same mesh cell.
Since $\hat \phi_h(\xi)$ is linear (for $\mathbb{P}_1$
elements) or quadratic (for $\mathbb{Q}_1$ elements)
along the line segment connecting the two points,
the root $\xi_\Gamma$ can be easily found using a closed-form
formula. The monotone sequence $\{\xi_i\}_{i=0}^m$
can also be constructed efficiently, especially if
the quadrature point
$\mathbf{x}_Q$ lies in a narrow band around $\Gamma_h(t)$.

\begin{rmk}
  The SDF property of $\phi_h$ is not required by the second
  version of our algorithm. The approximate normal $\bfq_h$ can
   be constructed by differentiating a volume-of-fluid 
  function $\phi_h^{\rm VOF}\approx H_\epsilon(\phi_h^{\rm LS})$ or
  a phase field indicator function $\phi_h^{\rm PF}\approx
  S_\epsilon(\phi_h^{\rm LS}$)
  instead of $\phi_h^{\rm LS}$ defined by
  \eqref{weakS} and \eqref{weakq}. 
  \end{rmk}

\subsection{Gradient reconstruction}
\label{sec:interp}

Having determined the location of $\mathbf{x}_\Gamma$, we
reconstruct the gradient at that point using  finite
difference approximations on a uniform stencil aligned
with the unit outward normal vector
$$\mathbf{n}_\Gamma=-\frac{\mathbf{q}_h(\mathbf{x}_{\Gamma})}{|\mathbf{q}_h(\mathbf{x}_{\Gamma})|}\approx\mathbf{n}_+(\mathbf{x}_{\Gamma}),$$
as shown in Fig.~\ref{sketch}(a). The
straight line $\hat{\mathbf{x}}(\xi)$ through $\mathbf{x}_Q
=\hat{\mathbf{x}}(0)$
intersects the inner
 edge of the diffuse interface region at the point
$\mathbf{x}_P=
\mathbf{x}_\Gamma-\epsilon\mathbf{n}_\Gamma$.
The three-point interpolation stencil
$$\mathcal S_{\rm 2D}(\mathbf{x}_\Gamma)=
\{\mathbf{x}_P-0.5\epsilon\boldsymbol{\tau}_\Gamma,\mathbf{x}_P,
\mathbf{x}_P+0.5\epsilon\boldsymbol{\tau}_\Gamma\}
$$
is constructed using a unit vector $\boldsymbol{\tau}_\Gamma$
orthogonal to $\mathbf{n}_\Gamma$. In the 3D case, the
five-point stencil
$$\mathcal S_{\rm 3D}(\mathbf{x}_\Gamma)=
\{\mathbf{x}_P-0.5\epsilon(\mathbf{n}_\Gamma\times\boldsymbol{\tau}_\Gamma),
\mathbf{x}_P-0.5\epsilon\boldsymbol{\tau}_\Gamma,\mathbf{x}_P,
\mathbf{x}_P+0.5\epsilon\boldsymbol{\tau}_\Gamma,
\mathbf{x}_P+0.5\epsilon(\mathbf{n}_\Gamma\times\boldsymbol{\tau}_\Gamma)
\}
$$
is used for interpolation purposes. Note that the cross product
$\mathbf{n}_\Gamma\times\boldsymbol{\tau}_\Gamma$ is a unit vector
orthogonal to $\mathbf{n}_\Gamma$ and $\boldsymbol{\tau}_\Gamma$.
Hence, the vectors $\boldsymbol{\tau}_\Gamma$ and
$\mathbf{n}_\Gamma\times\boldsymbol{\tau}_\Gamma$ span the
shifted tangential plane at $\mathbf{x}_P$.

 If a Neumann (rather than Dirichlet) boundary
condition is prescribed on $\Gamma(t)$, we use it to
define the interface value $\partial_n U(\mathbf{x}_\Gamma)$ directly.
Otherwise, we use the Dirichlet value $u_\Gamma(\mathbf{x}_\Gamma)$
and the interpolated value $u_h(\mathbf{x}_\Gamma-\epsilon\mathbf{n}_\Gamma)$ to construct the finite difference approximation
\beq\label{normalder}
\partial_nU(\mathbf{x}_\Gamma)=\frac{ u_\Gamma(\mathbf{x}_\Gamma)-
  u_h(\mathbf{x}_\Gamma-\epsilon\mathbf{n}_\Gamma)}{\epsilon}
\eeq
to the normal derivative at $\mathbf{x}_\Gamma$. Tangential
derivatives are approximated similarly using interpolated
values at the points belonging to the interpolation stencil
($\mathcal S_{\rm 2D}(\mathbf{x}_\Gamma)$ in 2D,
$\mathcal S_{\rm 3D}(\mathbf{x}_\Gamma)$ in 3D).
Since
the interface thickness parameter $\epsilon$ is
proportional to the mesh size $h$, the reconstructed
gradient $\nabla U(\mathbf{x}_\Gamma)$ represents a
consistent approximation to $\nabla u_h(\mathbf{x}_P)$
with $\mathbf{x}_P$ approaching $\mathbf{x}_\Gamma$
as $h\to 0$.

\blue{
\begin{rmk}
  If the boundary data function $u_\Gamma$ of the continuous problem
  is defined only on the zero level
  set $\Gamma(t)$ of the exact level set function $\phi$, it may not
  be possible to  evaluate  $u_\Gamma$
  at $\mathbf{x}_\Gamma$. However, in most
  practical applications $\phi$ is known analytically
  or a globally defined function $u_\Gamma:\Omega\to\R$ 
  has to be evaluated on the ``shifted boundary'' $\Gamma_h(t)$ 
  because the location of
  $\Gamma(t)$ is unknown.
\end{rmk}}
  
\subsection{Flux extrapolation}
\label{sec:fluxes}

To evaluate $G(\phi_h,u_h,u_\Gamma)$, as
defined by \eqref{Gflux}, at the point $\mathbf{x}_Q$,
we need the values of $\partial_nU(\mathbf{x}_Q)$ and
$U(\mathbf{x}_Q)$. Since we are using
a $\mathbb{P}_1$ or $\mathbb{Q}_1$ finite element approximation,
constant extrapolation \beq\label{extrapDu}
\partial_n U(\mathbf{x}_Q)=
\partial_n U(\mathbf{x}_\Gamma)\eeq
is sufficient for the normal derivative. The linear extrapolation
formula (cf. \cite{sbm1})
\beq\label{extrapU}
U(\mathbf{x}_Q)=u_\Gamma(\mathbf{x}_\Gamma)+\nabla U(\mathbf{x}_\Gamma)
\cdot(\mathbf{x}_Q-\mathbf{x}_\Gamma)
\eeq
\cyan{can be used for extension of interface values. The
 inviscid component $F$ of the flux $G(\phi_h,u_h,u_\Gamma)$
  can also be defined by linearly extrapolating
$f_n(u_h(\mathbf{x}_\Gamma),u_\Gamma(\mathbf{x}_\Gamma))$ as we do
  in Section \ref{sec:hyperbolic}.}

Note that extrapolation is performed along the vector
$\mathbf{x}_Q-\mathbf{x}_\Gamma$,
which is generally not parallel to the normal $\mathbf{n}_\Gamma$
(see Fig.~\ref{sketch}(b)). The gradient
$\nabla U(\mathbf{x}_\Gamma)$ to be used in \eqref{extrapU}
is uniquely defined by the reconstructed normal and
tangential derivatives. Its calculation requires a
coordinate transformation from the local coordinate
system aligned with $\mathbf{n}_\Gamma$ to the standard
Cartesian reference frame. The corresponding transformation
formulas can be found, e.g., in  \cite{vectorlim}.
No transformations need to be performed for the normal
derivatives  $\partial_n U(\mathbf{x}_\Gamma)$
to be used in \eqref{extrapDu}.

\subsection{Ghost penalties}
\label{sec:ghost}

The best way to define a ghost penalty function $u_{\Omega,h}$
for \eqref{weakh} depends on the problem at hand and on the available data.
A simple and natural Dirichlet extension is given by
\beq\label{extrapG}
u_{\Omega,h}(\mathbf{x}_Q)=\Heps(\phi_h(\mathbf{x}_Q))u_h(\mathbf{x}_Q)+
(1-\Heps(\phi_h(\mathbf{x}_Q)))U(\mathbf{x}_Q).
\eeq
Ghost penalties of Neumann type (as mentioned in Remark \ref{rem:ghost})
can be defined, e.g., using
\beq\label{extrapG2}
\mathbf{g}_{\Omega,h}(\mathbf{x}_Q)=\Heps(\phi_h(\mathbf{x}_Q))\nabla u_h(\mathbf{x}_Q)+
(1-\Heps(\phi_h(\mathbf{x}_Q)))\nabla U(\mathbf{x}_\Gamma).
\eeq
This definition corresponds to constant extrapolation of
$\nabla U(\mathbf{x}_\Gamma)$ into $\Omega_-(t)$.
The dimensions of local penalization parameters $\gamma^e=\gamma_{\Omega,h}|_{K^e}$
should be chosen in accordance with the type of the
extension (Dirichlet: $\gamma_D^e=\mathcal
O(h^{-1})$, Neumann:
$\gamma_N^e=\mathcal O(h)$). Detailed analysis of ghost
penalties is beyond the scope of the present work. However,
we envisage
that it can be performed following \cite{burman,olsh}.

Implicit treatment of ghost penalties changes the sparsity
pattern of the finite element matrix. For that reason, we implement
implicit  schemes\,/\,steady-state solvers
using fixed-point iterations
\begin{align*}
  L_h\big(u_h^{(k+1)}\big)+\int_{\Omega_h}\gamma_{\Omega,h}w_hu_h^{(k+1)}\dx
  &=R_h\big(u_h^{(k)}\big)+\int_{\Omega_h}\gamma_{\Omega,h}w_hu_{\Omega,h}^{(k)}\dx,\\
   L_h\big(u_h^{(k+1)}\big)+\int_{\Omega_h}\gamma_{\Omega,h}\nabla w_h\cdot\nabla u_h^{(k+1)}\dx
   &=R_h\big(u_h^{(k)}\big)+\int_{\Omega_h}\gamma_{\Omega,h}\nabla
   w_h\cdot\mathbf{g}_{\Omega,h}^{(k)}\dx
\end{align*}
based on a split form $L_h\big(u_h\big)=R_h\big(u_h\big)$ of the fully
discrete problem without penalization (using $\gamma_{\Omega,h}=0$).
One iteration may suffice if the problem is time dependent
and the time step is small.

\begin{rmk}\label{rmk-dgp-lump}
\cyan{In our numerical studies for the linear advection equation (see
  Section \ref{sec:hyperbolic}), we found it useful to perform mass
  lumping for $\int_{\Omega_h}\gamma_{\Omega,h}w_hu_h\dx$
  but not for $\int_{\Omega_h}\gamma_{\Omega,h}w_hu_{\Omega,h}\dx$.
  This treatment of ghost penalty terms
 makes the method more stable without adversely affecting its accuracy.}
  \end{rmk}
  
\subsection{Extension velocities}
\label{sec:extvel}

If the evolution of $\Gamma_h(t)$ is driven by interfacial phenomena
and only the normal velocity $v_n=V|_{\Gamma_h}$ is provided by the
mathematical model (as in \cite{hogea,vermolen2007}), a velocity
field $\mathbf{v}_h$ for evolving the level set function $\phi_h$
can be constructed using \eqref{Gflux} with a suitably defined
extension velocity $V$. The value
$V(\mathbf{x}_Q)$ can be determined using constant or linear 
extrapolation, as in \eqref{extrapDu} and \eqref{extrapU},
respectively. For example, suppose that $V(\mathbf{x}_\Gamma)
=\mu\partial_nC(\mathbf{x}_\Gamma)$ is proportional to a
(reconstructed) normal derivative $\partial_nC(\mathbf{x}_\Gamma)$
of a concentration field $C$. Then $V(\mathbf{x}_Q)=
V(\mathbf{x}_\Gamma)$ may be used to define 
$\mathbf{v}_h(\mathbf{x}_Q)=V(\mathbf{x}_Q)\mathbf{n}_Q$.
This way to calculate $\mathbf{v}_h(\mathbf{x}_Q)$ is
an efficient alternative to methods that perform a
constant extension of $V(\mathbf{x}_\Gamma)$ in the normal
direction
by solving a partial differential equation (cf. \cite{sethian,chen,utz}).

\subsection{Damping functions}
\label{sec:damping}

If the regularized delta function $\delta_\epsilon(\phi)$ has a global
support, extrapolation needs to be performed at each quadrature
point in the fictitious domain. In this case, the cost of numerical integration
and matrix assembly can be drastically reduced  using
damping functions such as
\beq\label{damp}
D_{\epsilon}(\phi)=H_{\epsilon}(\phi+m\epsilon)-
H_{\epsilon}(\phi-m\epsilon),\qquad m\ge 2.
\eeq
In particular, computations involving $G(\phi_h,u_h,u_\Gamma)$ can be restricted to
a narrow band region (cf. \cite{adalst}) if the compact-support extensions
$U_\epsilon=D_\epsilon(\phi_h)U$ and
$\partial_n U_\epsilon=D_\epsilon(\phi_h)\partial_nU$ are used in \eqref{Gflux}.

Localized  counterparts of the ghost penalty
extensions \eqref{extrapG} and \eqref{extrapG2}
can be defined as follows:
\begin{align*}
u_{\Omega,h}&=\Heps(\phi_h)u_h+
(1-\Heps(\phi_h))D_\epsilon(\phi_h)U,\\
\mathbf{g}_{\Omega,h}&=\Heps(\phi_h)\nabla u_h+
(1-\Heps(\phi_h))D_\epsilon(\phi_h)\nabla U.
\end{align*}
Note that the integrals of
$\gamma_{\Omega,h}w_h(u_h-u_{\Omega,h})$ and
$\gamma_{\Omega,h}\nabla w_h\cdot\nabla(u_h-u_{\Omega,h})$ vanish
also in subdomains where $\Heps(\phi_h)=1$. Thus an
efficient narrow-band implementation is possible.

\cyan{ For the 
  extension velocity $\mathbf{v}_h=(\phi_h) V\mathbf{q}_h$,
  damping is performed automatically in \eqref{weakSh}.
  Indeed, the integral of
$w_h\bfv_h\cdot\nabla\Seps(\phi_h)$
vanishes in elements in which
$\Seps(\phi_h)$ is constant. This compact support property of
the advective term speeds up computations.
Moreover, there is no
 need to prescribe unknown values of
$\phi$ on inflow boundaries of $\Omega$, while the
standard level set method requires so.}

\section{Test problems and results}
\label{sec:num}

In this section, we apply the diffuse level set method to three new test problems. We design them to have closed-form polynomial exact solutions, which are independent of time even if the equation is time-dependent and the interface is moving. Problems of elliptic, parabolic, and hyperbolic type are defined in Sections \ref{sec:elliptic}--\ref{sec:hyperbolic}. The numerical results presented in these sections were calculated using the finite element interpolant $\phi_h$ of the exact signed distance function $\phi$. In the elliptic case, the interface $\Gamma=\partial\Omega_+$ is a fixed embedded  boundary. The objective of our study is to quantify the error caused by weak imposition and the proposed numerical treatment of interface conditions. In particular, ghost penalties of Dirichlet and Neumann type, as well as the narrow-band damping functions \eqref{damp},
are tested in this experiment.  The parabolic and hyperbolic cases are designed to assess additional errors caused by the motion of the interface $\Gamma(t)$ and by convective fluxes across fixed embedded boundaries, respectively. In Section \ref{sec:full}, we solve the parabolic problem using \eqref{weakS} and \eqref{weakq} to evolve the level set function $\phi_h$ numerically. The normal velocity of the interface is extended into~$\Omega$ using extrapolation with damping (as  described in Sections \ref{sec:extvel} and \ref{sec:damping}). This experiment illustrates that all components of our diffuse level set algorithm fit together and produce reasonable results in 2D. 

Computations for all test cases are performed on uniform quadrilateral meshes using $\epsilon=2h$ and a  C\texttt{++} implementation of the methods under investigation in the open-source finite element library MFEM \cite{mfem}. The numerical solutions are visualized using the open-source software {\sc GlVis} \cite{glvis}. 

\subsection{Elliptic test}
\label{sec:elliptic}

The elliptic scenario corresponds to the steady-state limit of
\eqref{pde} with $\mathbf{f}(u) = 0$ and $\kappa = 1$. We consider
the fictitious domain formulation of the boundary value problem
\begin{alignat*}{2}
\Delta u &= 0&&\quad\text{in }\Omega_+,\\
u &= u_\Gamma &&\quad\mbox{on }\Gamma.
\end{alignat*}
The domain $\Omega_+$ is enclosed by
 $\Gamma = \{(x,y)\in \R^2: (x-0.5)^2 +
(y-0.5)^2 = 0.0625\}=\partial\Omega_+$ and embedded into 
$\Omega=(0,1)^2$. The boundary data and analytical solution
are given by
\begin{equation}\label{sol:ellipt}
u_\Gamma (x,y,t) = (x-0.5)^2 - (y-0.5)^2=u(x,y,t)\qquad\forall t\ge 0.
\end{equation}
In our numerical experiments,
we use the constant penalization parameter $\gamma_D=h^{-1}$
for the Dirichlet extension $u_{\Omega,h}$ and
 $\gamma_N=h$ for the Neumann extension
$\mathbf{g}_{\Omega,h}$. 

The $L^2$ error 
of the finite element approximation in the embedded domain
$\Omega_+$ is approximated by 
$$
\|u_h-u\|_{0,\Omega_+}=\sqrt{\int_{\Omega_h}\Heps(\phi)(u_h-u)^2\dx}.
$$
To keep the integration domain $\Omega_{+,\epsilon}
=\{\mathbf{x}\in\Omega\,:\,
\Heps(\phi(\mathbf{x}))>0\}$
constant on all refinement levels, we use the interface
thickness $\epsilon$ of the finest mesh in this formula.

\begin{figure}[h!]
\centering
\begin{minipage}[t]{0.5\textwidth}
  \centering (a)  $u_h$, damped DGP
\includegraphics[width=0.9\textwidth]{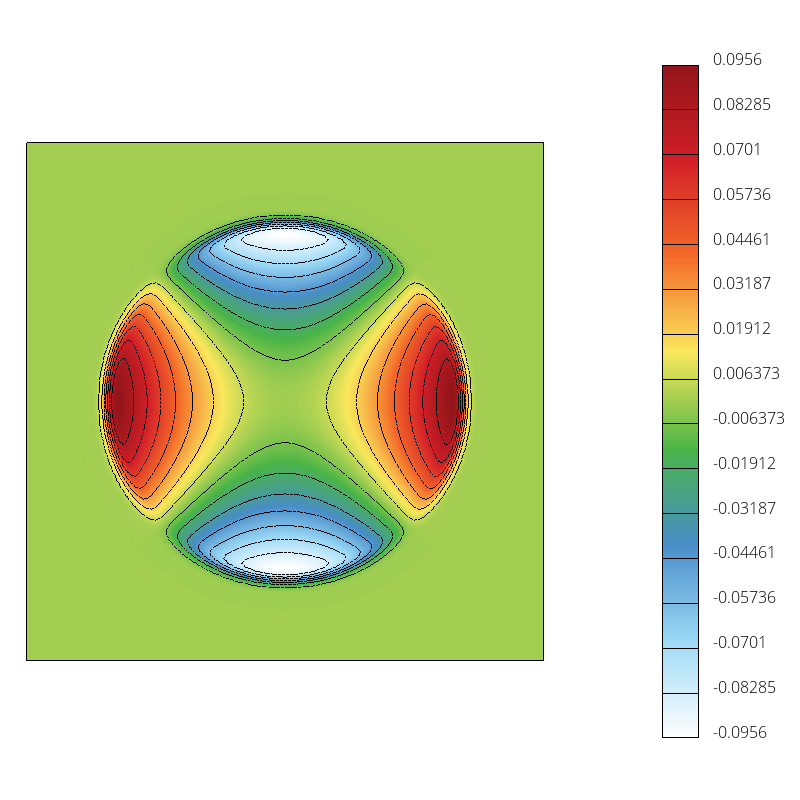}
\end{minipage}%
\begin{minipage}[t]{0.5\textwidth}
  \centering (b)  $H_\epsilon (\phi)u_h$, damped DGP
  \includegraphics[width=0.9\textwidth]{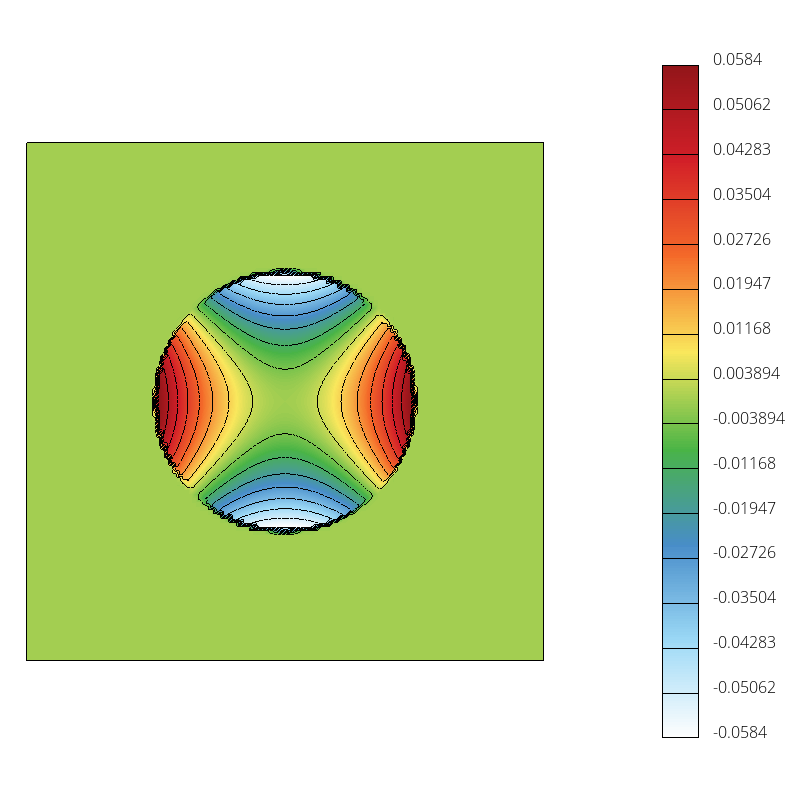}
\end{minipage}
\vskip0.5cm

\begin{minipage}[t]{0.5\textwidth}
  \centering (c) $u_h$, damped NGP
\includegraphics[width=0.9\textwidth]{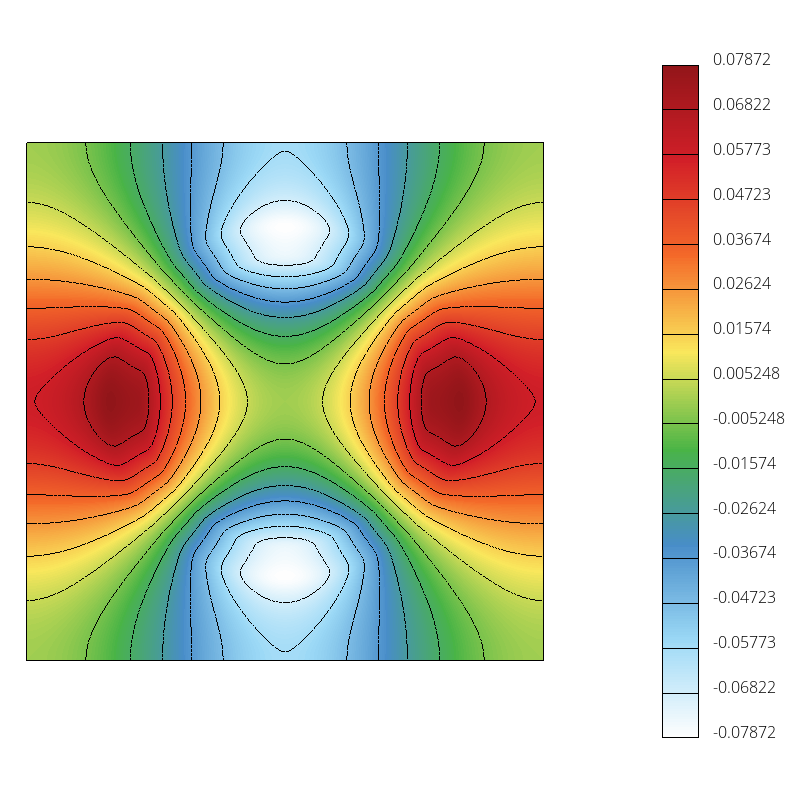}
\end{minipage}%
\begin{minipage}[t]{0.5\textwidth}
  \centering (d) $H_\epsilon (\phi)u_h$, damped NGP
  \includegraphics[width=0.9\textwidth]{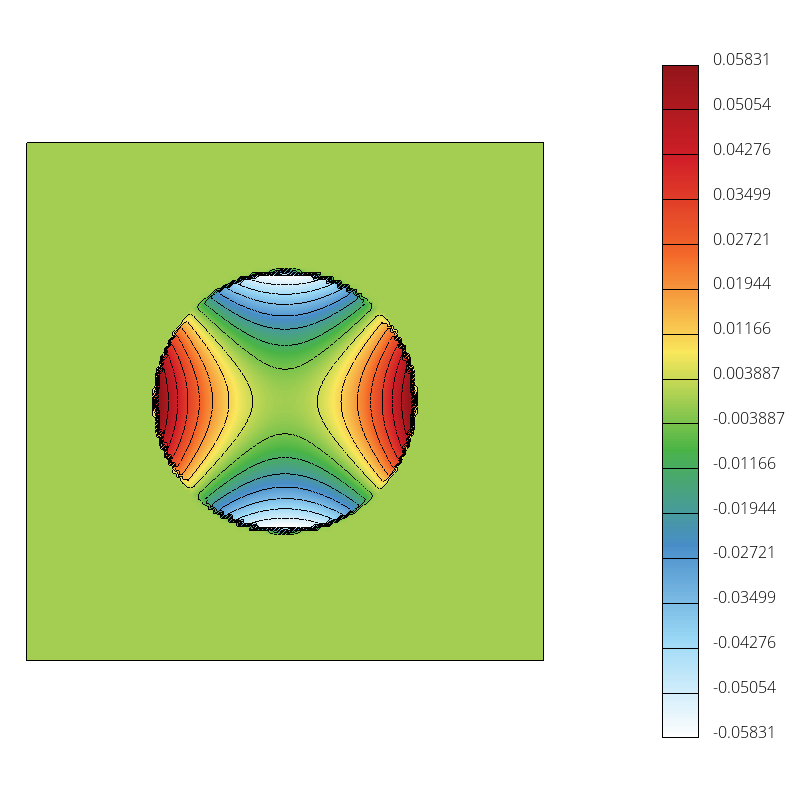}
\end{minipage}
\vskip0.25cm 
\caption{Elliptic test,
  numerical solutions obtained using $h=\frac{1}{128}$.}
\label{fig:ell}
\end{figure}

Table \ref{tab:EllEOC} summarizes the results of grid convergence
studies for the diffuse level set algorithm using Gaussian
regularization \eqref{delta:gauss} and four kinds of
ghost penalties. The
experimental order of convergence (EOC) is approximately 2
for all types of penalization that we compare in this study
(Dirichlet and Neumann, without and with multiplication
by the damping functions $D_{\epsilon}(\phi)$ defined by \eqref{damp}).
Hence, the restriction of level set extrapolation
to the narrow band $\Omega_\epsilon:=
\mathrm{supp}(D_{\epsilon}(\phi_h))$ around the embedded
boundary has no adverse effect
on the accuracy of the proposed method.
The numerical solutions obtained using damped ghost penalties
and the mesh size $h=\frac{1}{128}$ are shown in Fig.~\ref{fig:ell}.

\begin{table}[ht]
  \centering
  \begin{tabular}{ccccccccc}
  \hline
  $h^{-1}$ & full NGP & EOC & damped NGP & EOC & full DGP & EOC & damped DGP & EOC \\
  \hline 
  16   & 3.61e-03 &      & 3.61e-03 &      & 2.28e-03 &      & 2.28e-03 &      \\ 
  32   & 1.05e-03 & 1.78 & 1.06e-03 & 1.77 & 6.87e-04 & 1.73 & 6.87e-04 & 1.73 \\ 
  64   & 2.45e-04 & 2.10 & 2.61e-04 & 2.02 & 1.72e-04 & 2.00 & 1.72e-04 & 2.00 \\ 
  128  & 5.50e-05 & 2.15 & 6.03e-05 & 2.11 & 3.97e-05 & 2.11 & 3.97e-05 & 2.11 \\
  256  & 1.25e-05 & 2.14 & 1.39e-05 & 2.12 & 8.70e-06 & 2.19 & 8.70e-06 & 2.19 \\
  512  & 3.06e-06 & 2.03 & 3.42e-06 & 2.03 & 1.95e-06 & 2.15 & 1.95e-06 & 2.15 \\
  1024 & 7.37e-07 & 2.05 & 8.46e-07 & 2.01 & 4.21e-07 & 2.22 & 4.21e-07 & 2.22 \\
  \hline
  \end{tabular}
  \caption{Elliptic test,
    $L^2$ convergence history for Dirichlet ghost penalties (DGP) and Neumann ghost penalties (NGP). The full version extends boundary data into $\Omega$. The damped version
    extends into the narrow band $\Omega_\epsilon$.}
  \label{tab:EllEOC}
\end{table}

\subsection{Parabolic test}
\label{sec:parabolic}

The parabolic test problem is a time-dependent version of the
elliptic case. We solve
\eqref{pde} with $\mathbf{f}(u) = 0$ and $\kappa = 1$.
The initial-boundary value problem is formulated as follows:
\begin{alignat*}{2}
\pd{u}{t}-\Delta u &= 0&&\quad\text{in }\Omega_+(t),\\
u &= u_\Gamma &&\quad\mbox{on }\Gamma(t),\\
u &= u_0 &&\quad\mbox{in }\Omega_+(0).
\end{alignat*}
The moving boundary
$\Gamma(t) = \{(x,y)\in \R^2: (x-0.5)^2 + (y-0.5)^2 = (0.25+0.15t)^2\}$
is an expanding circle of radius
$r(t)=0.25+0.15t$ centered at $(x-0.5,y-0.5)$. The final time
$T=1$ is chosen to be short enough for $\Omega_+(T)$ to
stay embedded in $\Omega=(0,1)^2$.
As in the elliptic test, the exact solution $u(x,y)$ and the Dirichlet
boundary data $u_\Gamma (x,y)$ are given by formula \eqref{sol:ellipt}.

\begin{figure}[h!]
\centering
\begin{minipage}[t]{0.5\textwidth}
  \centering (a) $u_h$, damped DGP
\includegraphics[width=0.9\textwidth]{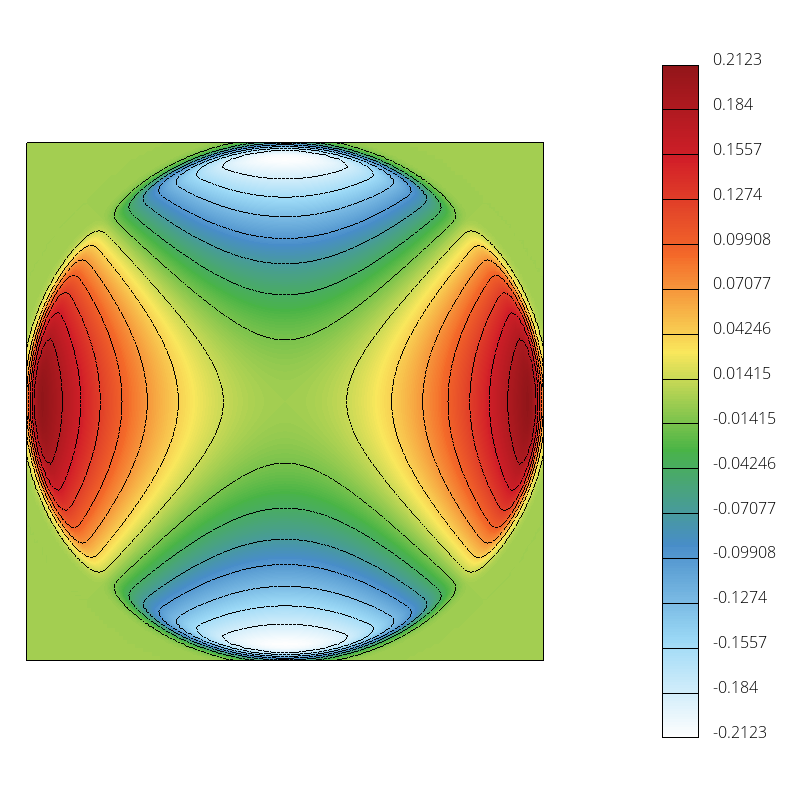}
\end{minipage}%
\begin{minipage}[t]{0.5\textwidth}
  \centering (b) $H_\epsilon (\phi)u_h$, damped DGP 
  \includegraphics[width=0.9\textwidth]{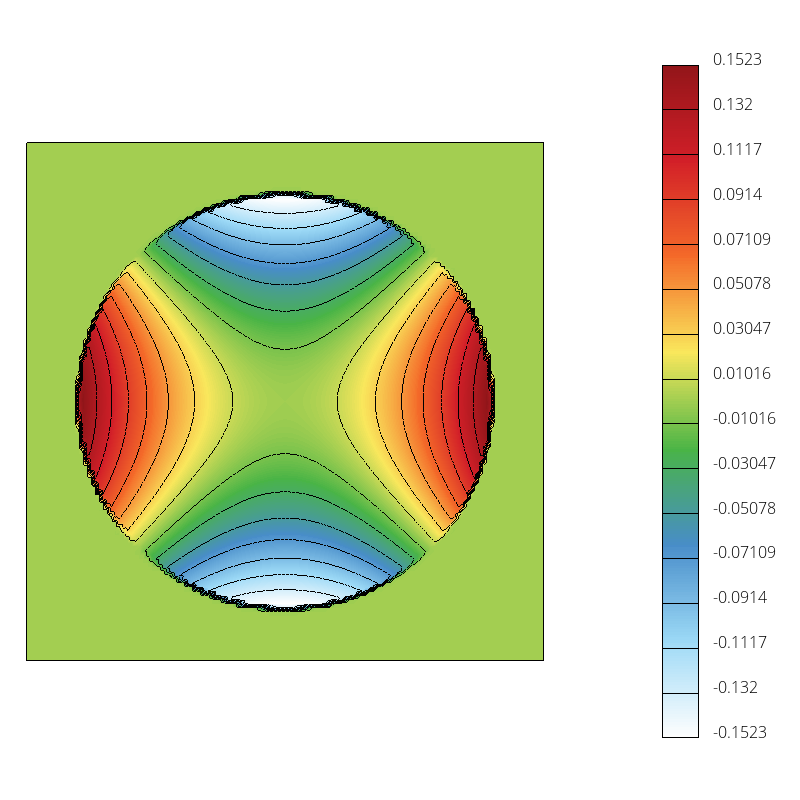}
\end{minipage}
\vskip0.25cm  
\caption{Parabolic test, numerical solution at $T=1$
  obtained using $h=\frac{1}{128}$.}
\label{fig:par}
\end{figure}

In this example, and in the following ones, we use Dirichlet ghost penalties
(DGP). The penalty parameter $\gamma_D$ is
defined as before. Time discretization is performed using the
second-order accurate Crank-Nicolson scheme and the time step
\cyan{$\Delta t=3.2h$} for the parabolic test.
The $L^2$ errors \blue{for numerical solutions
  at the final time and the corresponding}
EOCs are presented in Table \ref{tab:ParEOC}.
Second-order convergence behavior is observed again in this experiment.
The damped DGP result corresponding to the mesh size $h=\frac{1}{128}$ 
is shown in
Fig. \ref{fig:par}.

\begin{table}[ht]
  \centering
  \begin{tabular}{cccccc}
  \hline
  $h^{-1}$ & $(\Delta t)^{-1}$ & full DGP & EOC & damped DGP & EOC \\
  \hline 
  16   & 5   & 6.30e-03 &      & 6.30e-03 &      \\ 
  32   & 10  & 1.16e-03 & 2.44 & 1.16e-03 & 2.44 \\ 
  64   & 20  & 2.86e-04 & 2.02 & 2.86e-04 & 2.02 \\ 
  128  & 40  & 6.95e-05 & 2.04 & 6.95e-05 & 2.04 \\
  256  & 80  & 1.67e-05 & 2.06 & 1.67e-05 & 2.06 \\
  512  & 160 & 4.24e-06 & 1.98 & 4.24e-06 & 1.98 \\
  1024 & 320 & 1.10e-06 & 1.95 & 1.10e-06 & 1.95 \\
  \hline
  \end{tabular}
  \caption{Parabolic test,
    $L^2$ convergence history for Dirichlet ghost penalties (DGP). The full version extends boundary data into $\Omega$. The damped version
    extends into the narrow band $\Omega_\epsilon$.}
  \label{tab:ParEOC}
\end{table}

\subsection{Hyperbolic test}
\label{sec:hyperbolic}

The third test case is the hyperbolic ($\kappa=0$)  limit of
\eqref{pde} with the linear flux
$\mathbf{f}(u)=\mathbf{v}u$, where 
$$\mathbf{v}(x,y) = (0.5-y,x-0.5)^T$$
is a divergence-free rotating
velocity field. We solve the solid body rotation problem
\begin{alignat*}{2}
  \pd{u}{t} + \nabla\cdot(\mathbf{v}u)&= 0&&\quad\text{in }\Omega_+,\\
u &= u_\Gamma &&\quad\mbox{on }\Gamma_{\rm in},\\
u &= u_0 &&\quad\mbox{in }\Omega_+,
\end{alignat*}
where $\Gamma_{\rm in}$ is the inflow part of
$\Gamma(t) = \{(x,y)\in \R^2: (x - 0.75)^2 + (y - 0.5)^2 = 0.0225\}$
and $\Omega_+$ is
embedded into $\Omega=(0,1)^2$.
Computations are stopped at the
final time $T=1$. The exact solution 
\begin{equation}
u(x,y,t) = (x-0.5)^2+(y-0.5)^2\qquad \forall t\ge 0
\end{equation}
is used to define the boundary data $u_\Gamma(x,y,t)$ at the
inlet $\Gamma_{\rm in}$ when it comes to constructing $u_{\Omega,h}$.

\cyan{In this example,
the component $F$ of $G(\phi_h,u_h,u_\Gamma)$ is given by
the extended upwind approximation 
$$
  F_Q(t)=V_Q\hat U_Q,
    $$
    where $V_Q=\mathbf{v}(x_\Gamma,y_\Gamma)\cdot\mathbf{n}_\Gamma$
    and  $\hat U(x_Q,y_Q,t)$ is the result of linear extrapolation for
$$
\hat u(x_\Gamma,y_\Gamma,t)=\begin{cases}
 u_\Gamma(x_\Gamma,y_\Gamma,t) & \mbox{if}\
\mathbf{v}(x_\Gamma,y_\Gamma)\cdot\mathbf{n}_\Gamma<0,\\
 u_h(x_\Gamma,y_\Gamma,t) & \mbox{if}\ \mathbf{v}(x_\Gamma,y_\Gamma)\cdot\mathbf{n}_\Gamma\ge 0.
\end{cases}
$$}

\begin{figure}[h!]
\centering
\begin{minipage}[t]{0.5\textwidth}
  \centering (a) $u_h$, damped DGP
\includegraphics[width=0.9\textwidth]{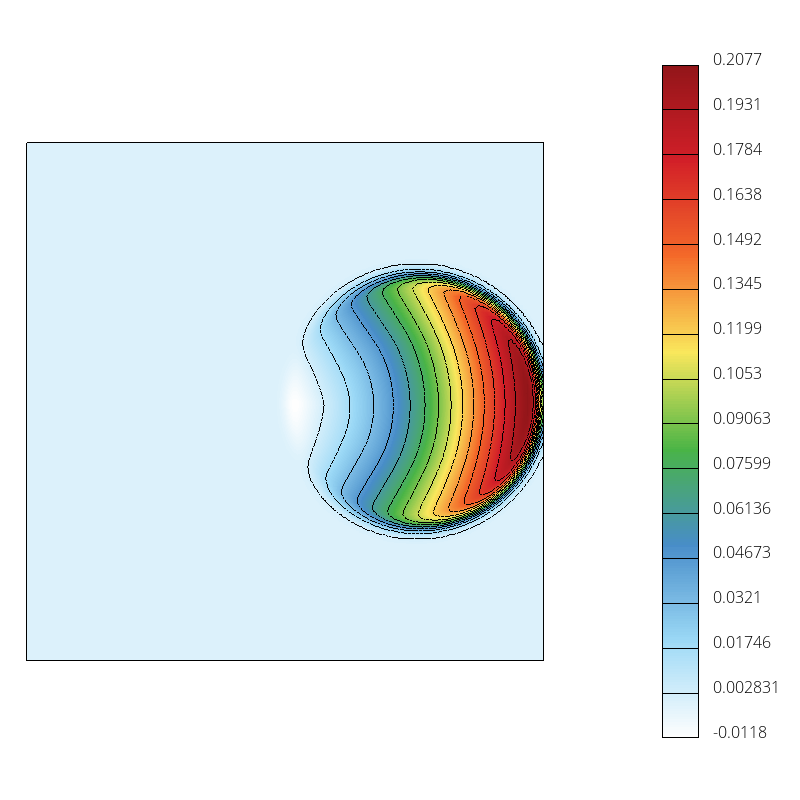}
\end{minipage}%
\begin{minipage}[t]{0.5\textwidth}
  \centering (b) $H_\epsilon (\phi)u_h$, damped DGP
  \includegraphics[width=0.9\textwidth]{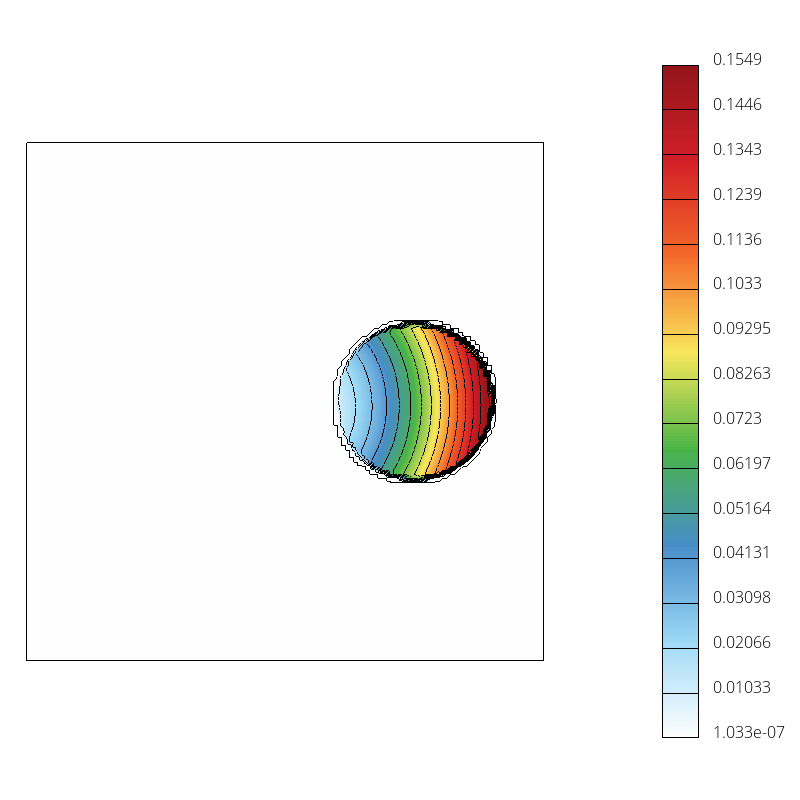}
\end{minipage}
\vskip0.25cm 
\caption{Hyperbolic test, numerical solution at $T=1$
  obtained using $h=\frac{1}{128}$.}
\label{fig:hyp}
\end{figure}

The parameter settings for the space discretization are chosen as
before. Time integration is performed using Heun's method, a
second-order explicit Runge-Kutta scheme which is known
 to be strong stability preserving (SSP) under a CFL-type
condition. We use the time step \cyan{$\Delta t=0.5h$} corresponding
to the maximum CFL number \cyan{$\nu=0.25$} in our computations.
The results are presented in Table \ref{tab:HypEOC}
and Fig. \ref{fig:hyp}. \cyan{In this
experiment, the consistent-mass version of the ghost
penalty term produced small ripples on fine
meshes. The use of mass lumping for $\int_{\Omega_h}
\gamma_{\Omega,h}w_hu_h\dx$ (as  mentioned in
Remark~\ref{rmk-dgp-lump}) cured this issue.
The first-order convergence behavior for
 the solid body rotation problem is not
surprising because we discretize a pure advection
equation using the continuous Galerkin method without any
stabilization other than mass lumping for the 
ghost penalty term.}

\begin{table}[ht]
  \centering
  \begin{tabular}{cccccc}
  \hline
  $h^{-1}$ & $(\Delta t)^{-1}$ & full DGP & EOC & damped DGP & EOC \\
  \hline 
  16   & 32   & 4.11e-03 &      & 4.11e-03 &      \\ 
  32   & 64   & 2.35e-03 & 0.81 & 2.35e-03 & 0.81 \\ 
  64   & 128  & 1.24e-03 & 0.92 & 1.24e-03 & 0.92 \\ 
  128  & 256  & 6.45e-04 & 0.95 & 6.45e-04 & 0.95 \\
  256  & 512  & 3.30e-04 & 0.97 & 3.30e-04 & 0.97 \\
  512  & 1024 & 1.67e-04 & 0.98 & 1.67e-04 & 0.98 \\
  1024 & 2048 & 8.40e-05 & 0.99 & 8.40e-05 & 0.99 \\
  \hline
  \end{tabular}
  \caption{Hyperbolic test,
    $L^2$ convergence history for Dirichlet ghost penalties (DGP). The full version extends
boundary data into $\Omega$. The damped version
    extends into the narrow band $\Omega_\epsilon$.}
  \label{tab:HypEOC}
\end{table}

\subsection{Level set advection}
\label{sec:full}

\begin{figure}[h!]
\centering
\begin{minipage}[t]{0.5\textwidth}
  \centering (a) constant extrapolation
  
\includegraphics[width=0.9\textwidth]{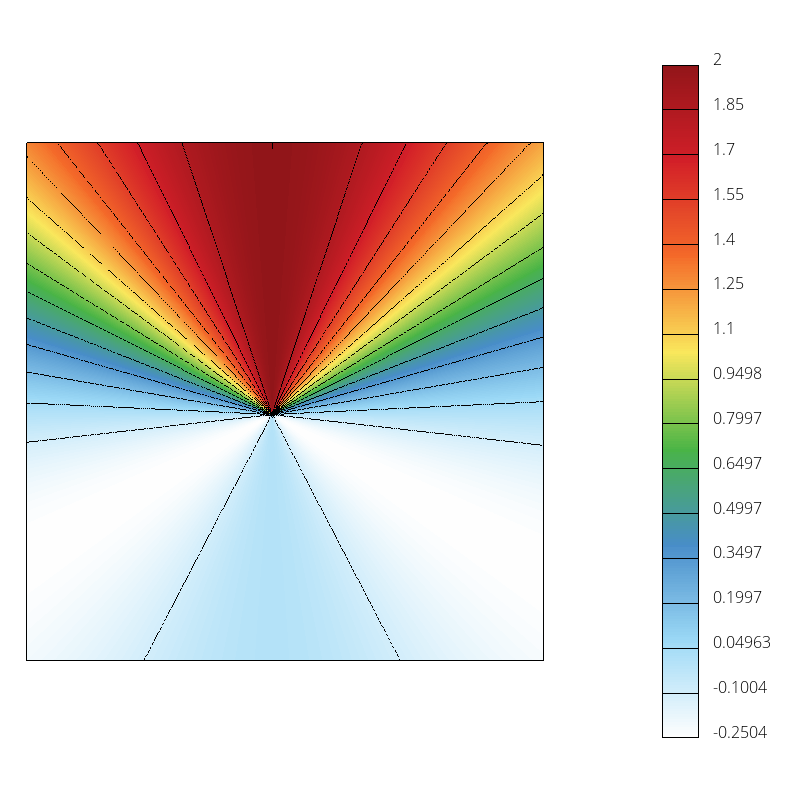}
\end{minipage}%
\begin{minipage}[t]{0.5\textwidth}
  \centering (b) elliptic extension
  
  \includegraphics[width=0.9\textwidth]{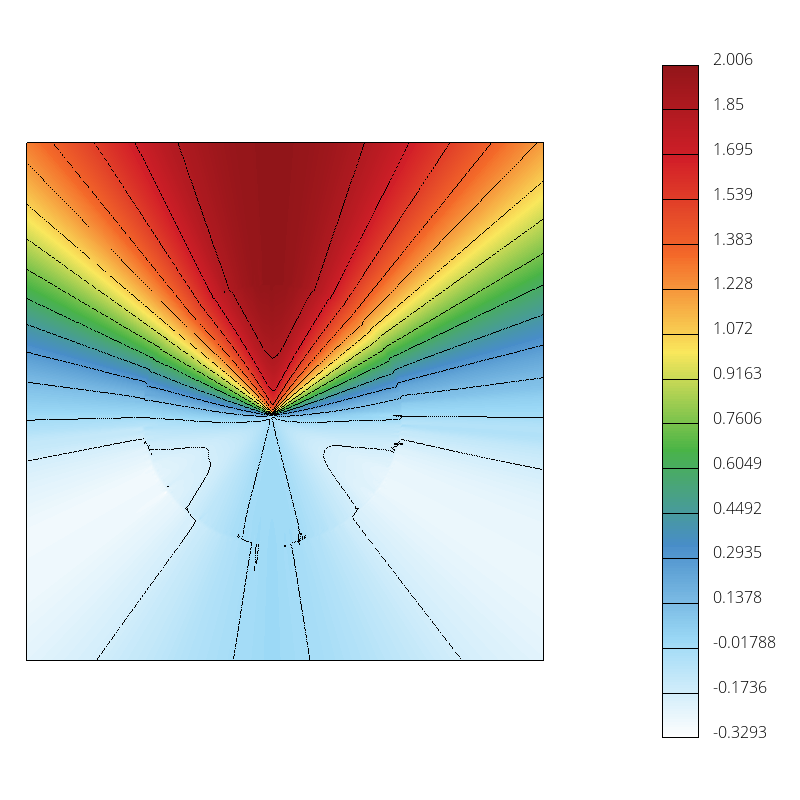}
\end{minipage}
\vskip0.25cm  
\caption{Circular test case \cite[Section 4.1]{utz},
  $h=\frac{1}{128}$, extensions
  constructed using (a) constant extrapolation, as proposed in
  Section~\ref{sec:extvel}, and (b)~a continuous
  Galerkin version of the elliptic extension method developed in \cite{utz}.}
\label{fig:ext}
\end{figure}

In the numerical studies above, we used the finite element interpolant
$\phi_h$ of the exact signed distance function $\phi$
for extrapolation purposes.
In the final test, we use the monolithic conservative
level set method \cite{mono} to compute
a numerical approximation $\phi_h$ for
the parabolic test problem formulated in
Section~\ref{sec:parabolic}.
The constant normal velocity $v_n=0.15$ of the interface
$\Gamma(t)$ is approximated by
$$
  V(x_\Gamma,y_\Gamma,t) = -0.15\partial_n\Phi(x_\Gamma,y_\Gamma\cyan{,t}),
$$
where $\partial_n\Phi(x_\Gamma,y_\Gamma\cyan{,t})$ is defined in terms of $\phi_h$
similarly to \eqref{normalder}. That is, we use
$$
\partial_n\Phi(x_\Gamma,y_\Gamma\cyan{,t})=\frac{0-
  \phi_h(x_P,y_P\cyan{,t})}{\epsilon}
=-\frac{1}{\epsilon}\phi_h(x_P,y_P\cyan{,t}).
$$
Since $\phi(\mathbf{x}_P\cyan{,t})=\epsilon$ for the exact SDF evaluated
at $\mathbf{x}_P=
\mathbf{x}_\Gamma-\epsilon\mathbf{n}_\Gamma$, the so-defined
approximate normal derivative should, indeed, satisfy
$\partial_n\Phi(x_\Gamma,y_\Gamma\cyan{,t})\approx -1$ for $\phi_h\approx\phi$.
We use constant extrapolation $V(x_Q,y_Q,t)=V(x_\Gamma,y_\Gamma,t)$
 to calculate the
extension velocities $$\mathbf{v}_h(x_Q,y_Q,t)=\cyan{-V(x_Q,y_Q,t)
\mathbf{q}_h(x_Q,y_Q,t)}
$$ at quadrature points in the narrow band. \cyan{Since
  multiplication by $\nabla\Seps(\phi_h)$ localizes
  the advective term in \eqref{weakSh}, no artificial
  damping functions are used for $\mathbf{v}_h$.}
To verify the accuracy
of the algorithm
proposed in Section~\ref{sec:extvel},
we applied it the circular test problem from
\cite[Section 4.1]{utz}. The exact constant extension
of $V_\Gamma(x,y)=y(1+y)$ defined on
$\Gamma=\{(x,y)\in\R^2\,:\,\sqrt{x^2+y^2}=1\}$
 is
$$
V(x,y)=y\frac{\sqrt{x^2+y^2}+y}{\sqrt{x^2+y^2}}.
$$
The results obtained with our method and with a
continuous Galerkin version of the elliptic extension procedure proposed
in \cite{utz} are shown in Fig.~\ref{fig:ext}. It can be seen that simple extrapolation performed remarkably well in this preliminary test.

\begin{figure}[h!]
  \centering\vspace{0.5cm}
  
\begin{minipage}[t]{0.5\textwidth}
  \centering (a) $H_\epsilon (\phi)u_h$
  
\includegraphics[width=0.9\textwidth]{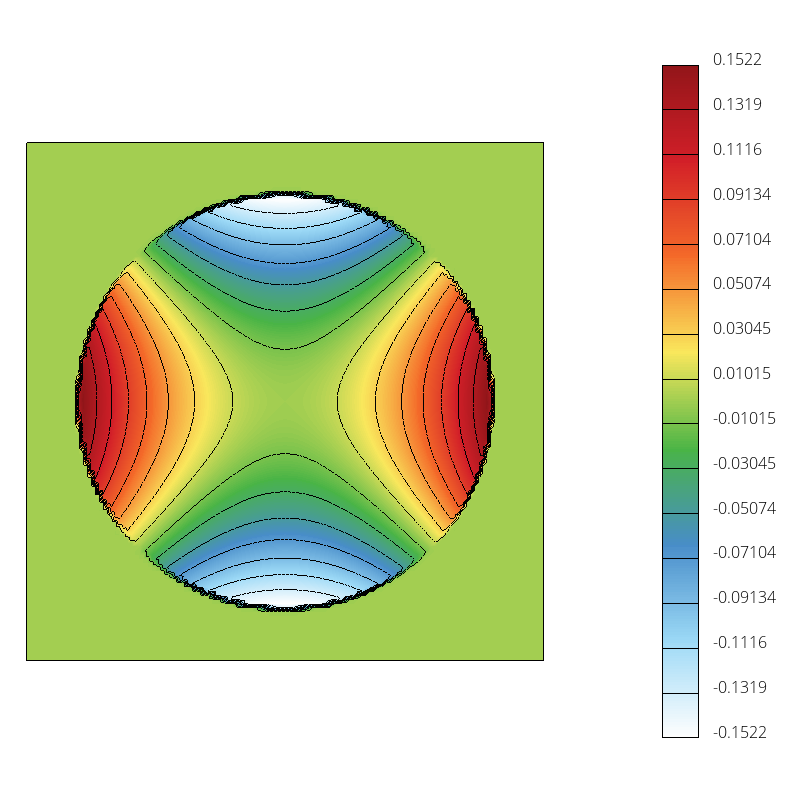}
\end{minipage}%
\begin{minipage}[t]{0.5\textwidth}
  \centering (b) $\phi_h$
  
  \includegraphics[width=0.9\textwidth]{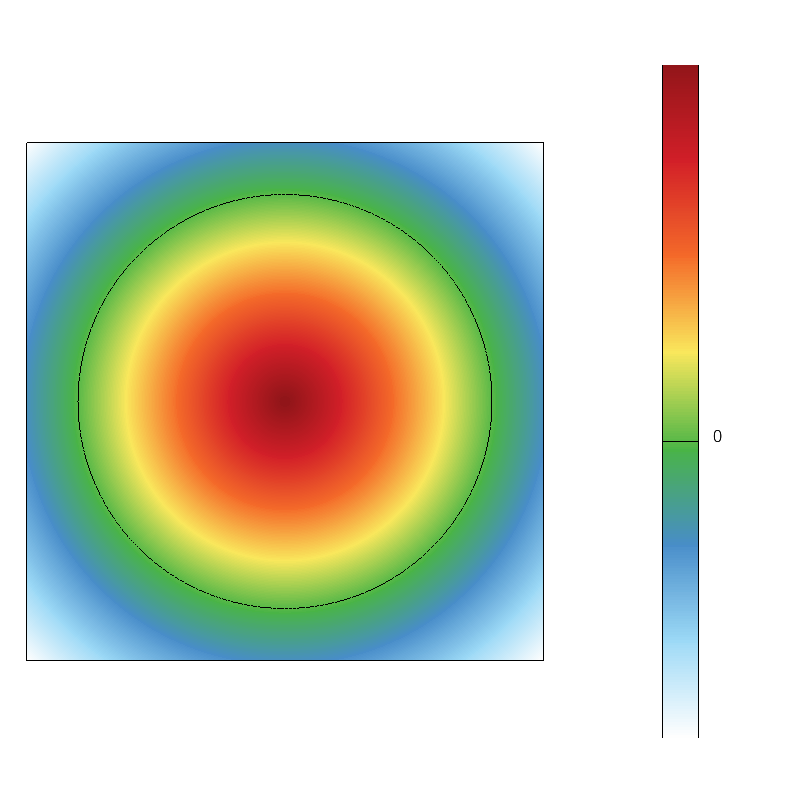}
\end{minipage}
\vskip0.25cm  
\caption{Parabolic test with numerical level set advection,
results obtained with $\epsilon=2h$ at $T=1$.}
\label{fig:lsa}
\end{figure}

\begin{figure}[h!]
  \centering\vspace{0.5cm}

  \includegraphics[width=0.9\textwidth]{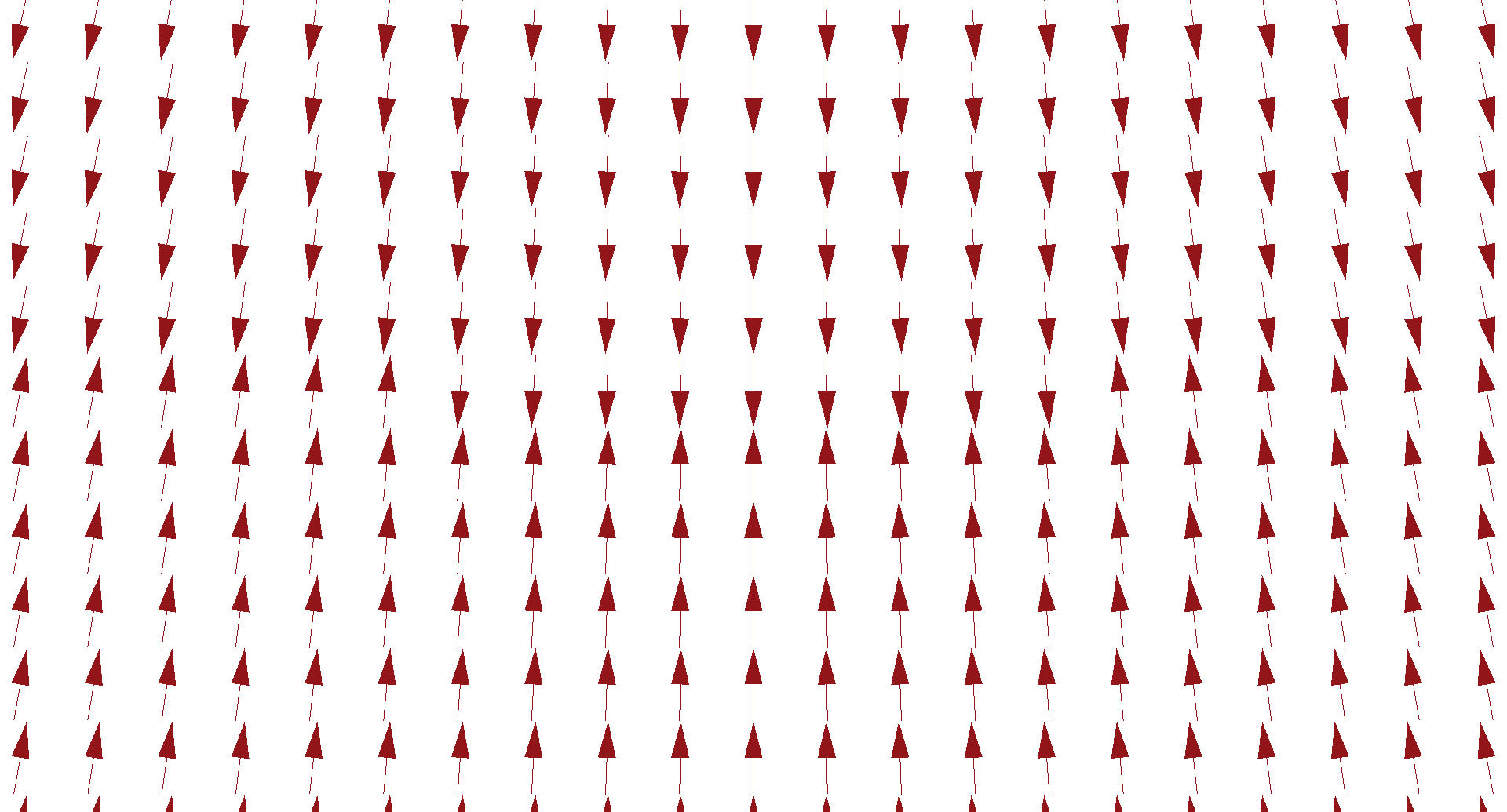}

\vskip0.25cm  
\caption{Nodal values of $\mathrm{sign}(\phi_h)\mathbf{q}_h$, a zoom for the bottom side of the interface at $T=1$.}
\label{fig:q_hplot}
\end{figure}

\begin{figure}[h!]
  \centering\vspace{0.5cm}
  
\begin{minipage}[t]{0.5\textwidth}
  \centering (a) $H_\epsilon (\phi)u_h$
  
\includegraphics[width=0.9\textwidth]{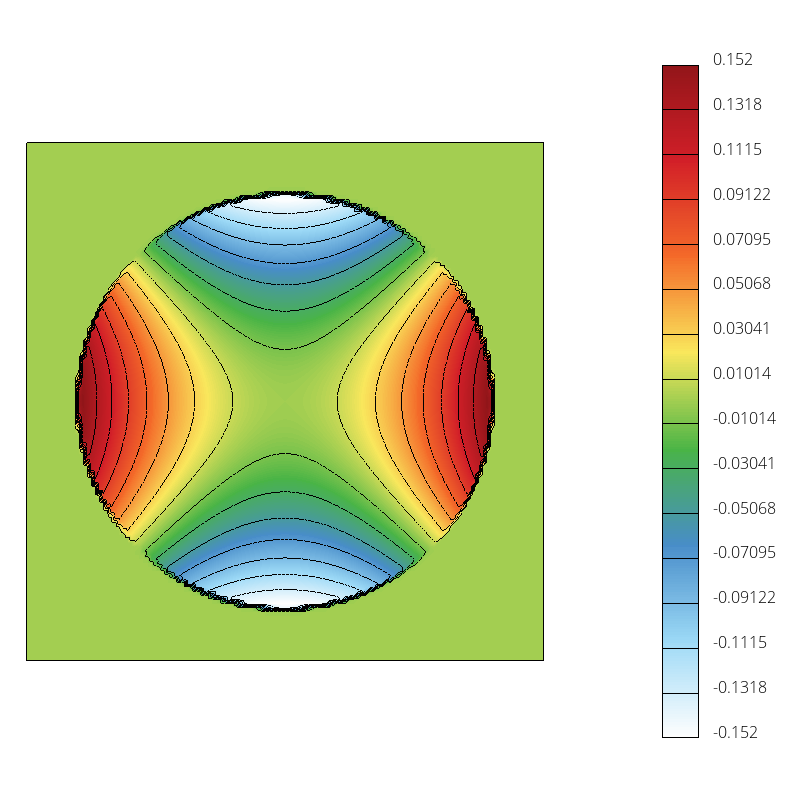}
\end{minipage}%
\begin{minipage}[t]{0.5\textwidth}
  \centering (b) $\phi_h$
  
  \includegraphics[width=0.9\textwidth]{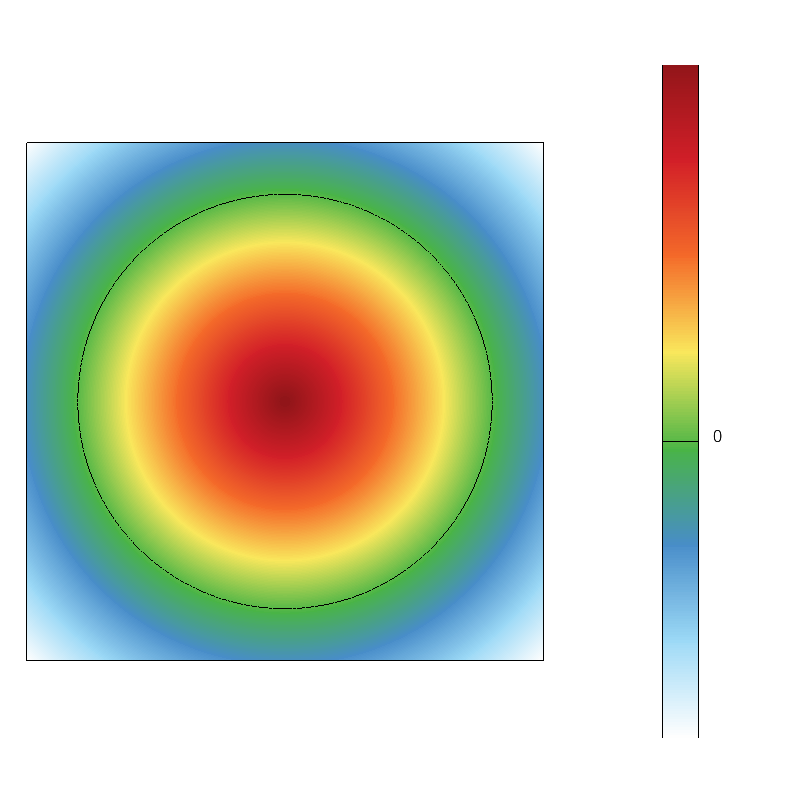}
\end{minipage}
\vskip0.25cm  
\caption{Parabolic test with numerical level set advection,
results obtained with $\epsilon=h$ at $T=1$.}
\label{fig:lsa2}
\end{figure}

In Fig.~\ref{fig:lsa} we present the numerical results for our
parabolic test with numerically advected level set function.
The snapshots show $u_h$ and \cyan{$\phi_h$} at
the final time. Both the structure of the exact solution
\eqref{sol:ellipt} and the
circularity of the interface are preserved reasonably well.
\blue{We also show the signed normal $\mathrm{sign}(\phi_h)\mathbf{q}_h$
  of the level set function in Fig.~\ref{fig:q_hplot}. It
  can be seen that all interface navigation vectors are
  indeed pointing towards $\Gamma(T)$.}
A direct comparison with Fig. \ref{fig:par} indicates that
disturbances caused by the error of numerical
approximation $\phi_h\approx\phi$ are of the same order
as perturbations caused by other numerical errors.
\blue{In our experience, the resuls are rather insensitive to the choice
  of the interface thickness parameter $\epsilon$
  (compare Fig.~\ref{fig:lsa} and Fig.~\ref{fig:lsa2}).}
We conclude
that our approach to construction of extension velocities is
a viable option for level set algorithms and other methods.

\section{Conclusions}
\label{sec:conclusions}

The unfitted finite element method presented in this paper
is based on a new fictitious domain formulation of an
initial-boundary value problem for a general time-dependent
conservation law. Sharp interface conditions are extended
into a diffuse interface. The main focus of our
investigations was on proper
definition and efficient numerical implementation of narrow-band
extension operators.
The general framework that we
developed for extrapolation of interfacial data features
\begin{itemize}
\item a fast closest-point search algorithm that uses a
 postprocessed gradient of a  level set function;
\item  a new way to define and calculate
 compact-support extensions of Dirichlet and Neumann
type;
\item narrow-band integration of terms involving extended fluxes,
  ghost penalties, and velocity fields.
\end{itemize}
The same methodology can be used for efficient extrapolation
to sharp surrogate boundaries (as in \cite{sbm1}). Theoretical
studies of related methods in \cite{burman,kublik2018,olsh,sbm,sbm1,zahedi}
and other publications provide useful tools for further analysis of
the proposed approach. Other promising research avenues include
construction of new approximate delta functions with compact supports,
extensions to interface problems with jump conditions (as in \cite{sbm,wadbro}),
\blue{special treatment of surfaces with kinks and corners (as in
\cite{kublik2018}),} as well as theory-based design of ghost penalties
using Dirichlet and Neumann extensions.

\section*{Acknowledgments}

The first author would like to thank Prof. Guglielmo Scovazzi
(Duke University) for a personal introduction to the shifted
boundary method and inspiring discussions of the proposed
approach.

This work was supported by the German Research
Association (DFG) under grant KU 1530/28-1.

\end{document}